%% file: main.tex
\documentclass[letterpaper,onecolumn,  journal]{IEEEtran}  
\pdfoutput=1                                                        

\IEEEoverridecommandlockouts                              

\usepackage{microtype}
\usepackage{graphicx}
\usepackage{subfigure}
\usepackage{booktabs} 
\usepackage{tabularx}
\usepackage{amsmath,amsfonts,bm}
\usepackage{amsthm, amssymb, bbm}
\usepackage{thmtools}
\usepackage{thm-restate}
\usepackage{lipsum}
\usepackage{mathtools}

\newtheorem{lemma}{Lemma}

\newtheorem{theorem}{Theorem}

\usepackage{hyperref}

\usepackage{graphics} 
\usepackage{graphicx}
\usepackage{array}
\usepackage{epsfig} 
\usepackage{color}
\usepackage{float}
\usepackage[noabbrev]{cleveref}
\usepackage{algorithm}
\usepackage{algorithmic}
\usepackage{url}
\usepackage{calrsfs}
\usepackage{textcomp}
\usepackage{multirow}
\usepackage{booktabs}

\usepackage{multirow}
\usepackage{hhline}
\DeclareMathAlphabet{\pazocal}{OMS}{zplm}{m}{n}
\DeclareMathAlphabet{\mathpzc}{OT1}{pzc}{m}{it}
\newcolumntype{P}[1]{>{\centering\arraybackslash}p{#1}}
\newcolumntype{N}{@{}m{0pt}@{}}

\newcommand{\bd}{\mathbf}
\newcommand{\R}{\mathbb{R}}

\begin{document}

	\title{Learning to Solve Network Flow Problems via Neural Decoding}
	
	\author{Yize Chen and Baosen Zhang \vspace{-4ex}
		
		\thanks{Y. Chen  and B. Zhang are with the Department of Electrical and Computer Engineering at the University of Washington,  emails: \{yizechen, zhangbao\}@uw.edu. }
	}
	
	\maketitle

\begin{abstract}
	\input{abstract}
\end{abstract}

\section{Introduction}
\label{sec:intro}
\input{intro}
\section{Modeling and Learning Formulation}
\label{sec:setup}
\input{setup}

\section{Learning Active Constraints }
\label{sec:kkt}
\input{kkt_derivation}

\section{Training Surrogate Neural Networks}
\label{sec:training}
\input{training}

\section{Robustness to Errors}
\label{sec:learning}
\input{Learning}

\section{Experimental Results}
\label{sec:experiment}
\input{simulation}

\section{Conclusion and Discussion}
\label{sec:conclusion}
\input{Conclusion}

\bibliographystyle{IEEEtran}
\bibliography{ref}

\newpage

\appendix
\input{appendix}

\end{document}

%% file: abstract.tex
Many decision-making problems in transportation, power system and operations research require repeatedly solving large-scale linear programming problems with a large number of different inputs. For example, in energy systems with high levels of uncertain renewable resources, tens of thousands of scenarios may need to be solved every few minutes. Standard iterative algorithms for linear network flow problems, even though highly efficient, become a bottleneck in these applications. In this work, we propose a novel learning approach to accelerate the solving process. By leveraging the rich theory and economic interpretations of linear programming duality, we interpret the output of the neural network as a noisy codeword, where the codebook is given by the optimization problem's KKT conditions. We propose a feedforward decoding strategy that finds the optimal set of active constraints, and the optimal solution can then be found by solving a linear equation. This design is error correcting and can offer significant speedup compared to current state-of-the-art iterative solvers, while providing much better solutions in terms of feasibility and optimality compared to other learning approaches.


%% file: intro.tex
In many engineering applications, optimization programs are solved repeatedly to make real-time decisions. Among them, network flow problems form an important class and have been studied for decades with wide-ranging applications. They arise naturally in the context of transportation, networking, communication and energy systems. In many of these settings, a network flow problem takes the form of a linear program (LP)~\cite{bertsimas1997introduction}.

Despite decades of studies, LPs can still face computational challenges in some settings. Often these challenges are due to the increased stochasticity and the low-latency requirements of real-time applications. The motivating application of this paper is in power systems, where intermittent and random renewable resources are being integrated into the grid at all levels~\cite{wen2015optimal, dalal2016hierarchical}. These new resources can create new flow patterns that fundamentally transform how systems operate. Inadequate planning for sudden events can lead to significant losses of welfare, as demonstrated by the recent rolling blackouts in Texas~\cite{Jimenez19} and California~\cite{Penn19}.

The resource allocation process in power systems is called the optimal power flow (OPF) problem, where it takes the form of a network flow problem that minimizes the cost of power generation to satisfy the loads, subject to all of the physical network constraints (e.g., generator limits and line capacities)~\cite{glover2012power}. Typically, a linear version called DCOPF--an LP problem--is often applied in practice and solved periodically (e.g., every 5 minutes) to find the optimal operating conditions~\cite{stott2009dc}. Because of the uncertainties brought by the renewables on many of the nodes, the number of generation and load scenarios that need to be considered are starting to grow exponentially. 
Even if each scenario under consideration takes $0.1$ seconds to solve using modern solvers, not all of them can be completed within the required time period. Therefore, using neural networks to learn the mapping between input load profiles and the corresponding optimal generation outputs has gained significant attention, since making inference via a trained architecture can be potentially orders of magnitude faster than an iterative solver~\cite{deng2016probabilistic,pan2019deepopf}.

\begin{figure*}[t]
	\centering
	\includegraphics[width=0.99\linewidth]{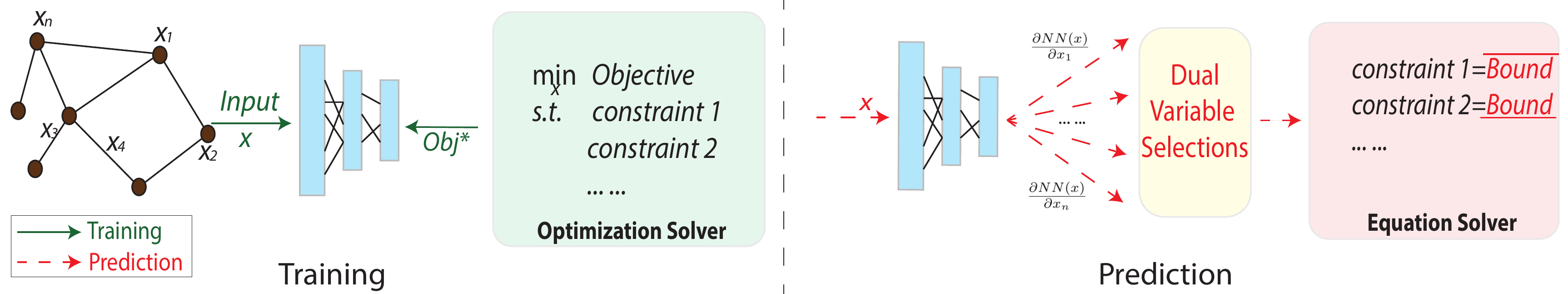}
	\caption{\footnotesize The schematic of our proposed Neural Decoder for solving LP problems. A neural network is trained to predict the optimal objective value; during implementation for solving LPs, the network's gradient is interpreted as a noisy codeword for active constraints, and a linear equation is solved to obtain optimal solutions. }
	\label{fig:workflow}
\end{figure*}

The learning-theoretic question becomes whether neural networks (NNs) can be used to directly make decisions in optimization problems, and in particular, \emph{whether it can learn the solutions of LPs as a function of the changes in the problem data}. Surprisingly, the answer to this question has been largely negative. In~\cite{amos2017optnet,agrawal2019differentiable}, convex optimization problems are embedded as layers to standard neural networks, with the rationale that directly learning the mapping from problem data to solutions is difficult, especially when learned solutions need to satisfy all of the constraints.

In this paper, we answer the question in the affirmative for LPs.
Instead of viewing it as an end-to-end learning task, we leverage the rich algorithmic understanding of LPs as well as the economic interpretation of the primal and dual variables to offer a novel solution architecture. Our workflow is shown in Figure~\ref{fig:workflow}. Concretely, we construct a neural network that takes the net load at each node as the input and outputs the optimal system cost (a scalar). Of course, the optimal value of the cost is not the solution of the LP. Rather, using the neural network, we compute the \emph{gradient} of the cost with respect to the net loads. Identifying these as the \emph{dual variables}, we use them to predict the binding nodal and line constraints. Once these constraints are identified, the optimal solution is given by solving a simple linear system of equations. The overall procedure can be seen as an efficient surrogate learning model for optimization solvers.

We make the following contributions in this work:
\begin{enumerate}
  \item We provide a novel architecture to solve network flow problems. By using a neural network to learn the mapping between the problem data and the optimal cost and \emph{interpreting the gradient of the network as dual variables}, we provide extremely efficient ways to identify the active constraints in the problem.
  \item Our method can be thought as decoding a system where the codebook is given by the KKT conditions of the LP. It is \emph{error correcting}, since even if output of the neural network is noisy, the correct set of active constraints at optimal solutions can still be identified.
  \item We show that our method can provide more than an order of magnitude speedup compared to current state-of-the-art iterative solvers and much better performance in terms of feasibility and optimality compared to other learning approaches.
\end{enumerate}

\subsection{Related Works and Background}

\noindent \textbf{Solving LP in Engineering Applications} \quad Generic LP solvers have been optimized to the point that they provide the state-of-the-art algorithms for many applications. For example, power system operators use CPLEX or Gurobi for their dispatch problems~\cite{bertsimas2012adaptive,bixby2015computational,Heidel16,gomez2018electric}. Specialized problems such as max flow have dedicated algorithms that have faster theoretical guarantees~\cite{dantzig2003max}, but general network flow problems will often use generic solvers~\cite{bertsimas1997introduction,vanderbei2015linear}. A key challenge is that many applications now require problem to be solved repeatedly with stringent latency constraints~\cite{rao2019engineering}.

\noindent \textbf{Learning to Solve Optimization} \quad  Our work falls under the category of using machine learning to solve optimization problems, which has a long history dating back to at least~\cite{kennedy1988neural, lillo1993solving}. A line of research considers using deep learning to solve combinatorial optimization problems typically by developing new heuristics~\cite{khalil2017learning, gasse2019exact}. This paper is closer to the work that tries to solve convex optimization problem by directly using supervised learning to find the optimal mapping from input data to the optimal solution, with applications in power systems~\cite{pan2019deepopf, canyasse2017supervised}, scheduling~\cite{cao2018machine, cui2019spatial}, and resource management~\cite{sun2017learning}. These algorithms can be regarded as end-to-end behavior cloning of expert policies~\cite{bojarski2016end}.

The challenge of using end-to-end models is that they can suffer from compounding errors and poor generalization performances, especially when hard constraints need to be satisfied. In~\cite{amos2017optnet}, a small LP problem with 3 variables is shown to be hard to learn via supervised learning. Our work is also related to \cite{deka2019learning, ng2018statistical}, where authors proposed to predict the set of active constraints at optimality.
The learned neural networks in these settings have difficulties generalizing because of the high dimensionality of their outputs and need training data to essentially cover input space.

In \cite{dalal2016hierarchical, mao2016resource, song2018learning}, the authors describe reinforcement learning formulation to model the interaction between learners' decisions and rewards, but there is no guarantee the solutions satisfy all constraints. In \cite{donti2017task, agrawal2019differentiable, djolonga2017differentiable, lee2019meta}, differentiable convex optimization layers provide convenient venues of designing end-to-end models, yet they still use standard solvers at runtime and do not achieve our goal of speeding up computations in time/resource-constrained applications.

%% file: setup.tex
In this section, we provide the formulation of two network flow problems and describe the relationship between the optimal solutions and the active constraints.

\vspace{-10pt}
\subsection{Network with Independent Edge Flows}
\vspace{-5pt}
A common type of network flow problems arising in logistics and communication network concerns the distribution of goods from origins~(e.g., a manufacturing plant) to destinations (e.g., consumers)~\cite{bertsimas1997introduction,srikant2013communication}. Given a graph $G=(V,E)$, the sources, destinations, and intermediate points are collectively modeled as nodes in the set $V=\{1,...,n\}$. The sources have certain amounts of goods and each of the sink has some demand that needs to be satisfied.
The goods need not be sent directly from source to destination and may be routed through intermediary points. For each node, let $x_i$ denote its production and let $l_i$ denote the load. We allow a node to have both positive $x_i$ and $l_i$ (a node can both produce and consume goods), so $\bd x=(x_1,\dots,x_n)$ and $\boldsymbol \ell=(l_1,\dots,l_n)$ are in $\mathbb{R}^n$.
 Each of the edge in the network carries some flow of goods and they are related to the nodes through a conservation equation: the sum of the flows into a node must equal to its net load. Algebraically, suppose there are $m$ edges, each with flow $f_j$. Then the flows $\bd f=(f_1,\dots,f_m)$ are related linearly to $\bd x$ and $\boldsymbol\ell$ through an incidence matrix $\bd A$ of the graph $G$, where $\bd x+ \bd A \bd f =\boldsymbol\ell$ and we follow the convention that flows into a node is positive.

We associate a positive $c_i$ with the $i$'th node, interpreted as the unit cost of producing the good. The standard network flow problem is then to minimize the total costs of production:
\begin{subequations}
	\label{eqn:network}
	\begin{align}\vspace{-10pt}
	\min_{\bd x, \bd f} \quad & \bd c^T \bd x\\
	\text { s.t. } \quad		&\bd 0 \leq \bd x \leq \bar{\bd x} \label{eqn:network_x}\\
 & -\underline{\bd f} \leq \bd f \leq \bar{\bd f} \label{eqn:network_f}\\
	&\bd x+\bd A \bd f= \boldsymbol\ell \label{subequ:network2}
	\end{align}
\end{subequations}
where $\bd x$ and $\bd f$ are the optimization variables with upper and lower bounds as shown in \eqref{eqn:network_x} and \eqref{eqn:network_f}, respectively. Note, if a node is a pure sink, we can set $\bar{x}_i=0$. 


  \begin{figure}[t]
	\begin{minipage}[c]{0.55\textwidth}
		\centering
		\includegraphics[width=0.9\linewidth]{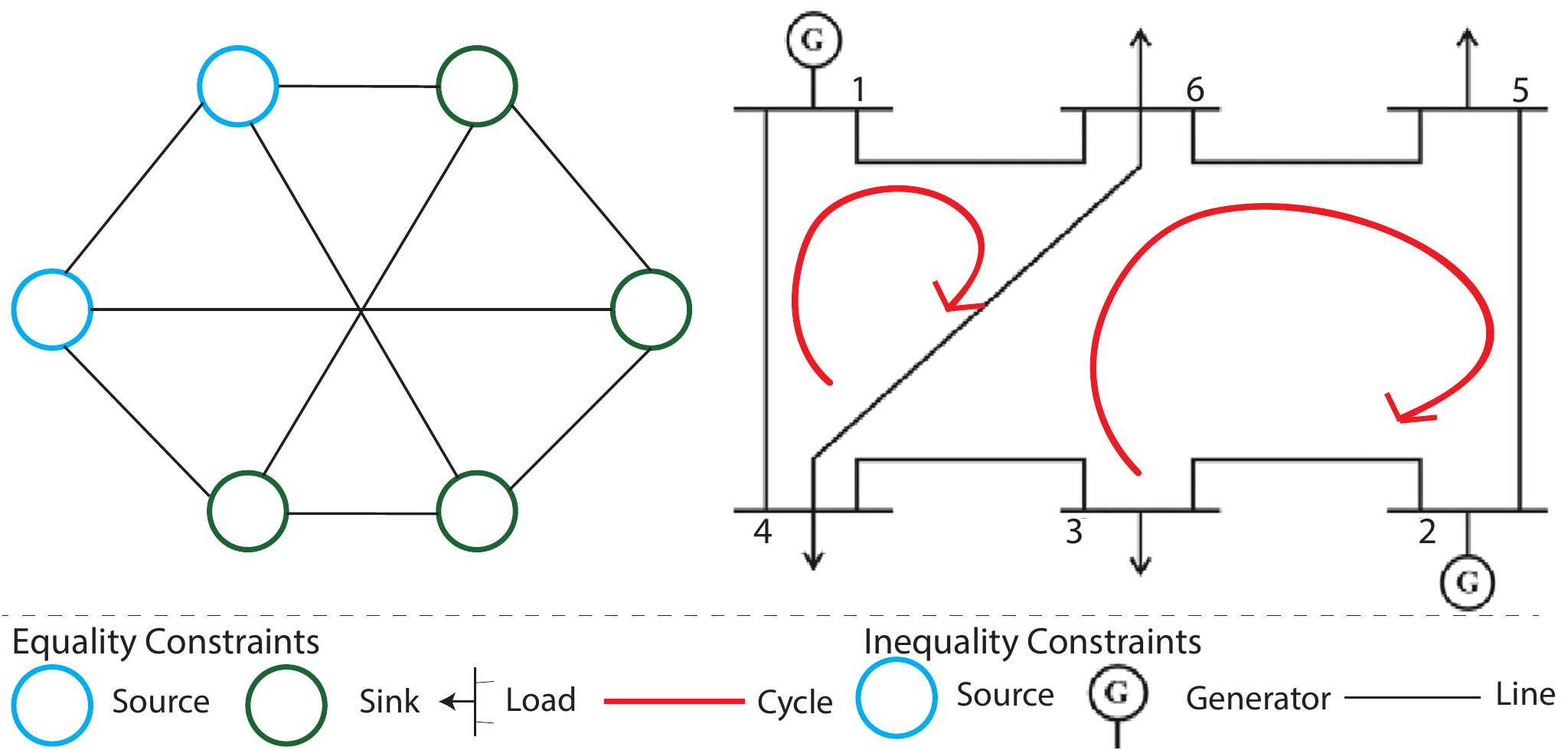}
	\end{minipage}\hfill
	\begin{minipage}[c]{0.4\textwidth}
		\caption{The network flow model and OPF model  along with physical constraints. Note that cycle equality constraints exist in OPF model due to the underlying power flow equations.}
		\label{fig:model}
	\end{minipage}
  \vspace{-10pt}
\end{figure}
\vspace{-10pt}
\subsection{Optimal Power Flow Problem}
\vspace{-5pt}
A network flow problem that includes constraints induced by physical laws is the optimal power flow (OPF) problem. It is one of the most important problems in power system operations, which finds the generations that minimizing cost while satisfying all the loads and flow constraints. The key difference with the problem in \eqref{eqn:network} is that the flows in a power system cannot be set independently. Namely, they must obey \emph{Kirchoff's voltage laws}, which states that a linear combination of the flows must sum up to zero in a cycle in the network as shown in Fig.~\ref{fig:model}. More modeling details can be found in Appendix~\ref{app:cycle}. Because of these \emph{cycle constraints}, the edge flows must lie in a subspace in $\mathbb {R}^m$. This is the space of fundamental flows, with a dimension of $n-1$. The rest of the $m-n+1$ flows are uniquely determined by the cycle constraints.  By selecting a spanning tree of the network where flows on the edges in the tree are fundamental, the rest of the flows on the edges are determined through the cycle constraints~\cite{gross2005graph}. Algebraically, we use $\bd f$ to denote the flows in a cycle basis and a matrix $\bd K \in \R^{m \times (n-1)}$ to map it to all of the flows. Then the optimization problem becomes:
\begin{subequations}
	\label{eqn:opf_DC}
	\begin{align}\vspace{-13pt}
	\min_{\bd x, \bd f} \quad & \bd c^T \bd x\\
	\text { s.t. } \quad & \bd 0\leq \bd x \leq \bar{\bd x} \label{eqn:opf_DC_gen}\\
	& -\underline{\bd f} \leq \bd K \bd f \leq \bar{\bd f} \label{eqn:opf_DC_flow}\\
	&\bd x+\bd{\tilde{A}} \bd f= \boldsymbol\ell \label{eqn:opf_DC_equ},\vspace{-10pt}
	\end{align}
\end{subequations}
where $\bd{\tilde{A}}$ is the modified incidence matrix mapping the fundamental flows to the nodal injections.



Both \eqref{eqn:network} are \eqref{eqn:opf_DC} are LPs and an instance of them can be solved efficiently using a solver like CVX or CPLEX. However, as it becomes necessary to resolve them for a large set of scenarios (e.g., to accommodate renewable energy), we need to explore a faster alternative. We do this by identifying the set of constraints that are active, or binding, at the optimal solution of the LPs.

\subsection{Active Constraints}
At the optimal solution of an LP, a constraint is called \emph{active} if it is an equality constraint or an inequality constraint that binds. In this paper, we assume the LP is not degenerative, and the key result (Theorem~\ref{thm:active}) we use from linear programming is that the number of active constraints is exactly the same as the number of decision variables in the problem.
\begin{theorem}
		\label{thm:active}
	For a nondegenerative linear programming problem with $p$ optimization variables, the optimal solution is determined by $p$ linearly independent active constraints.
\end{theorem}
This theorem states that if we know which constraints are active at the optimal, we can simply solve linear equations to find the optimal solution satisfying all constraints instead of solving an optimization problem. This is especially advantageous for network flow problems, where the number of constraints are typically much larger than the number of optimization variables.
The proof of the theorem is standard and can be found in many textbooks on linear programming~\cite{bertsimas1997introduction}.

In this paper, we develop a novel machine learning algorithm to learn the set of active constraints without solving the LP. Notably, we do not attempt to directly learn the set of active constraints via a classification problem since there can be exponential number of combination of active constraints, this approach has proven to be difficult~\cite{ng2018statistical}. Alternatively, we view active constraints as codes in a codebook, and the the learning process as decoding a noisy observation of a codeword.

%% file: kkt_derivation.tex
In this section, we present learning the active constraints as a decoding process, where the code can be learned implicitly.  By interpreting the KKT conditions as a codebook, we demonstrate how to utilize the structure and parameter information of the optimization model to train a neural network. Analysis in this section can also be extended to more general convex optimization problems.

Throughout this section, we work with the optimal power flow problem in \eqref{eqn:opf_DC}, since the standard network flow problem in \eqref{eqn:network} can be thought as a special case by taking $\bd K$ to be identity and the cycle basis to include flows on all of the edges. Case of line flow costs can be analyzed similarly.

\subsection{Derivatives of the LP}
We can view the LP in \eqref{eqn:opf_DC} as a mapping from the load vectors $\boldsymbol\ell$ to the optimal value of the objective function. We denote this function as $J^*(\boldsymbol\ell)$, which implicitly encodes the feasibility and optimality conditions. Since the output of $J$ is a scalar and is continuous in $\boldsymbol\ell$, it is much easier to learn than the vector of optimal solutions or the active constraints. Of course, we are rarely satisfied about knowing only the optimal cost. But this function has a rich structure that can be exploited using the following theorem (proof given in Appendix \ref{proof:mu}):
\begin{theorem}
	\label{thm:mu}
	Let $\bm \mu^* \in \mathbb{R}^n$ denote the optimal dual variables associated with equality constraints \eqref{eqn:opf_DC_equ}. If for a given $\boldsymbol\ell$  \eqref{eqn:opf_DC} is feasible, then 	$\nabla_{\boldsymbol\ell} J^*=\bm \mu^*$.
	Furthermore, the optimal solution for \eqref{eqn:opf_DC} is associated with the following active/inactive constraints:
	\begin{equation}
	\label{eqn:OPF_Constraints}
	x_i=\left\{
	\begin{array}{@{}ll@{}}
	\bar{x}_i, & \text{if }  \mu^*_i-c_i>0 \\
	0 & \text{if }  \mu^*_i-c_i<0\\
	\left(0,\,\bar{x}_i\right), &\text{otherwise}
	\end{array}\right.
	\end{equation}
\end{theorem}

Theorem \ref{thm:mu} states if $\bm \mu^*$ is known, then all active constraints at each of the node can be readily ``decoded'' by exploiting the parameters of original optimization problem. Of course, two questions remain: 1) how do we find the active \emph{flow constraints} and 2) how do we find $\bm \mu^*$?

\subsection{Finding Active Flow Constraints}
We now describe how to find the active constraints of the edge flows in \eqref{eqn:opf_DC_flow} which can be generalized to cases in  \eqref{eqn:network_f}. Let $\bar{\bm \nu}$ and $\underline{\bm \nu}$ be the dual variables associated with the capacity constraints in \eqref{eqn:opf_DC_gen}, $\bar{\bm \lambda}$ and $\underline{\bm \lambda}$ are the dual variables associated flow capacities in~\eqref{eqn:opf_DC_flow}, and $\bm \mu$ is the dual variable associated with the equality constraint~\eqref{eqn:opf_DC_equ}. The dual of \eqref{eqn:opf_DC} is:
\begin{subequations}
	\label{eqn:OPF_dual}
	\begin{align}
	\max_{\bm \mu, \bar{\bm \lambda},\underline{\bm \lambda}, \bar{\bm  \nu},\underline{\bm \nu}} \quad  &
	\mathbf{\bm \mu}^T\ \boldsymbol\ell
	-\underline{\bm \lambda}^T\underline{\bd f}-\bar{ \bm \lambda}^T\bar{\bd f}
	-\bar{\bm \nu}^T\bar{\bd x} \label{eqn:OPF_dual_1}\\
	\text { s.t. } \quad &\bd c-\bm \mu-\underline{\bm \nu}+\bar{\bm \nu}=\bd 0\label{eqn:OPF_dual_2}\\
	&-\tilde{\bd A}^T\bm \mu-\bd K^T\underline{\bm \lambda}+\bd K^T\bar{\bm \lambda}=\bd 0 \label{eqn:OPF_dual_3}\\
	&\bar{ \bm \nu}\geq \bd 0,\underline{\bm \nu}\geq \bd 0,
	\bar{\bm \lambda}\geq \bm 0,\underline{\bm \lambda}\geq  \bm 0
	\end{align}
\end{subequations}
If $\bm \mu^*$ is given (we will talk about how to learn this in Section~\ref{sec:training}), the dual variables associated with the nodal constraints and the edge flow constraints decouple, where the latter is:
\begin{subequations}
	\label{eqn:OPF_line_dual}
	\begin{align}
	\max_{\bar{{\bm \lambda}},\underline{\bm \lambda}} \quad  &-\underline{\bm \lambda}^T\underline{\bd f}-\bar{\bm \lambda}^T\bar{\bd f} \label{eqn:OPF_line_dual_1}\\
	\text { s.t. } \;  &  \bd K^T(\bar{\bm \lambda}-\underline{\bm \lambda})=\tilde{\bd A}^T\bm \mu^*  \label{eqn:OPF_line_dual_3}\\
	& \bar{\bm \lambda}\geq \bm 0,\, \underline{\bm \lambda}\geq  \bm 0
	\end{align}
\end{subequations}

We cannot directly find line constraint binding conditions due to the coupling between $\bar{\bm \lambda}, \underline{\bm \lambda}$ and $\bm \mu^*$.
The following simple lemma is a useful way to rewrite \eqref{eqn:OPF_line_dual}.
\begin{lemma}\label{lem:upsilon}
	Assume that $\underline{\bm f}=\bar{\bm f}$ (symmetric edge capacities). Then the problem \eqref{eqn:OPF_line_dual} is equivalent to
	\begin{subequations}
		\label{eqn:OPF_line_L1}
		\begin{align}
		\min_{\bm \upsilon} \quad  & ||\bm \upsilon||_1 \label{eqn:OPF_line_dual_L1}\\
		\text { s.t. } \;  &  \hat{\bd K}^T(\bm \upsilon)=\tilde{\bd A}^T\bm \mu^*   \label{eqn:OPF_line_L1l_3}
		\end{align}
	\end{subequations}
	where $\upsilon_i=\bar{f}_i \cdot \bar{\lambda_i}-\bar{f}_i \cdot \underline{\lambda_i}$, $|\upsilon_i|= |\bar{f}_i \cdot \bar{\lambda_i}-\bar{f}_i \cdot \underline{\lambda_i}|$, and $\hat{\bd K}=\mbox{diag}(1/\bar{f}_1,\dots,1/\bar{f}_m) \bd K$.
\end{lemma}
For the proof see Appendix~\ref{app:upsilon}. The assumption that the edge capacity is symmetric holds for most undirected networks seen in practice.

For the standard network flow problem in \eqref{eqn:network}, the matrix $\bd K$ is the identity and \eqref{eqn:OPF_line_L1l_3} is a full rank linear equation. Therefore, the optimal $\bm \upsilon$ (and hence the multipliers) can be found by inspection. However, for the OPF problem in \eqref{eqn:opf_DC}, the cycle constraints make solving \eqref{eqn:OPF_line_L1} less trivial.
Fortunately, transformed $\mathpzc{L}_1$ minimization problem in \eqref{eqn:OPF_line_L1} is extremely well studied, since it falls exactly into the regime of sparse signal recovery in compressed sensing, where signal $\bm \upsilon$ needs to be recovered via observation $\tilde{\bd A}^T\bm \mu^*$. Note that in sparse recovery, $\mathpzc{L}_1$ minimization is a surrogate problem where the ultimate goal is to find the sparsest solution to the problem. For us, \eqref{eqn:OPF_line_L1} is the exact problem we want to solve to finish the active constraints identification.

Because there are limited number of active constraints which is equivalent to limited number of nonzero entries in $\bm \upsilon$, while $\hat{\bd K}$ is very sparse, there are many algorithms that can solve \eqref{eqn:OPF_line_L1} extremely efficiently. For example, a family of greedy algorithms such as iterative hard thresholding~(IHT) can be used to find $\bm \upsilon$~\cite{blumensath2008iterative} in a few (fixed) number of iterations.

In our decoding scheme, KKT conditions and dual problems are first utilized to decode a set of dual variables as the status of inequality constraints, while the remaining constraints can be identified by a smaller, subsequent $\mathpzc{L}_1$ minimization step.






%% file: training.tex
By following Theorem~\ref{thm:mu}, we now describe the learning procedure of $J^*$ along with optimal dual variable $\bm \mu^*$.
The training is implemented offline using supervised data by collecting solution data of LP with different $\boldsymbol\ell$. We use the optimal cost, optimal multiplier (readily available from most solvers) and the active constraints set as the labeled data.

 We do not directly use a neural network to learn $\bm \mu^*$. Since $\bm \mu^*$ is not continuous in $\boldsymbol\ell$, learning $\bm \mu^*$ directly becomes a large classification problem, which is difficult because of the large number of possible values. Rather, we fit a neural network (parameterized by $\theta$) $g_\theta( \boldsymbol\ell)$ that maps optimization model input $ \boldsymbol\ell$ to the optimal cost. This function is continuous and piecewise linear, with distinct "breakpoints", where the derivative changes value. Therefore, it is naturally parameterized by using ReLU activation units. Let $h(g_\theta(\boldsymbol\ell))$ denote the binary vector indicating the active constraints determined by learner's solution based on the procedure described in the last section. The structure and parameters of the underlying optimization problem provide us with a number of terms in the design of learning objectives:
\begin{enumerate}[noitemsep,nolistsep]
	\item Regression loss defined between  $g_\theta(\boldsymbol\ell)$ and $J^*(\boldsymbol\ell)$ over the neural network's output;
	\item Regression loss for optimal dual variable defined between  $\nabla_{\boldsymbol\ell} g_\theta(\boldsymbol\ell)$ and $\bm \mu^*$;
	\item Distance loss defined between $h(g_\theta(\boldsymbol\ell))$ and $s^*(\boldsymbol\ell)$, where $s^*(\cdot)$ is the binary vector indicating active constraints at the optimal solution.
\end{enumerate}
The training loss is the sum of these terms
\begin{equation}
\label{eqn:loss}
\mathpzc{L}(\theta)= \| g_\theta(\boldsymbol\ell)- J^*(\boldsymbol\ell)  \|_2^2 +
\gamma_1 \| \nabla_{\boldsymbol\ell} g_\theta(\boldsymbol\ell) - \bm  \mu^* \|_2^2 +\gamma_2 \| h(g_\theta(\boldsymbol\ell))  -s^*(\boldsymbol\ell) \|_H
\end{equation}
where $\| h(g_\theta( \boldsymbol\ell))  -s^*(\boldsymbol\ell) \|_H $ is  the Hamming distance between active constraint sets\footnote{For binary vector $\mathbf{v}$, Hamming norm $\| \mathbf{v}\|_H$ is defined as the number of non-zero entries of vector $\mathbf{v}$}~\cite{norouzi2012hamming}, and $\gamma_1, \gamma_2$ are penalty parameters.

%% file: Learning.tex
The discussion in Section \ref{sec:kkt} motivates us to use a neural network to learn the optimal costs along with a ``codeword" associated with dual variables. Because we are dealing with continuous values, the prediction of the neural network will invariably have errors. Therefore, it is important that the errors made do not add up and cause incorrect identification of the active constraints. Using an analogy from communication theory, we think of the derivatives of the learned neural network $g_{\bm \theta}(\boldsymbol\ell)$ as \emph{noisy versions of a codeword}, and provide an  error-correcting approach to decode active constraints that is robust to errors made by the neural network. The overall algorithm is presented in Algorithm~\ref{alg:algorithm}.

Observe that for a given $\bm \mu^*$ and some noise $\bm \delta$, though the optimal solution for the optimization problem \eqref{eqn:OPF_line_dual} may be different for $\bm \mu^*$ and $ \bm \mu^* +\bm \delta$, the set of active constraints at optimal solutions can remain the same. Therefore, there exists a region around a ground truth optimal dual solution where as long as the noise does not push the solution outside of this region, the set of active constraints remains the same. It turns out these regions are polytopes and easily characterized by solving another linear program. This falls under the well studied area of linear programming sensitivity analysis~\cite{bertsimas1997introduction, jansen1997sensitivity}. The next lemma helps us to utilize these results to obtain an estimate of the testing error based on the training error of the neural network:
\begin{lemma}
	\label{lem:error}
	Consider a given $\boldsymbol\ell$ and its associated optimal dual variables $\bm \mu^*$. There exists a polytope $\mathpzc{P}_{\bm \mu^*}$ around $\bm \mu^*$ such that if $\tilde{\bm \mu} - \bm \mu^* \in \mathpzc{P}_{\bm \mu^*}$, then the active constraints determined by the procedure in Algorithm \ref{alg:algorithm} given $\tilde{\bm \mu}$ is correct. The set $\mathpzc{P}_{\bm \mu^*}$ is computable by a linear program.
\end{lemma}
The proof of this lemma is given in Appendix~\ref{app:error}. It shows that solutions via our approach is robust to errors, since we do not require that the value of predictions $\nabla_{\boldsymbol\ell}g_\theta( \boldsymbol\ell)$ to be exact, and there is no compounding of errors as compared to other end-to-end learning approaches. It also provides a way to check whether training is good enough: we can compare the error between the learned $\bm \mu$ and the actual $\bm \mu^*$, and if it falls within the polytope $\mathpzc{P}_{\bm \mu^*}$, then we are confident that the training results would be accurate. This lemma also leads to a way to use dictionary learning to further speed up the solution process, which we present in Appendix~\ref{app:dictionary}.

\begin{figure*}[hbt!]
	\centering
	\includegraphics[width=1.0\linewidth]{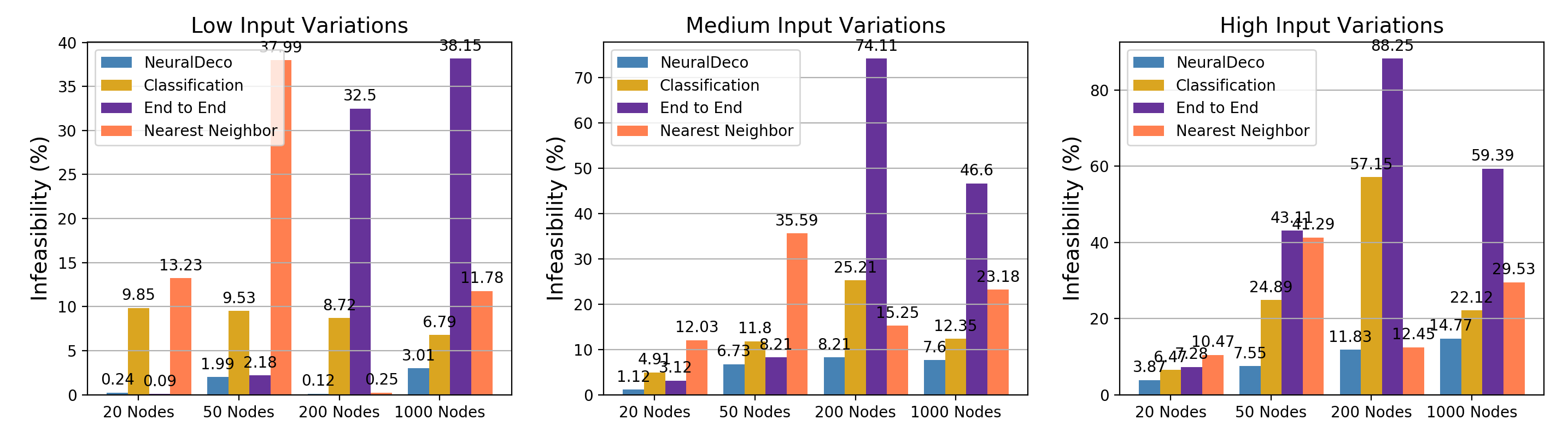}
	\caption{\footnotesize  The optimal active line constraints prediction results for Network Flow problem on different size of graphs. For varying problem input distributions (low, medium and high variances), our proposed Neural Decoder approach can learn the optimal set of active constraints with lowest infeasibility percentage comparing to active set classification, end-to-end prediction of line flow values and nearest neighbor prediction of active constraints.}
	\label{fig:NF}
\end{figure*}

%% file: simulation.tex

We evaluate the proposed learning approach (Neural Decoder) on a set of network flow optimization tasks, and demonstrate improvements over other methods. Specifically, over a wide range of problem input settings, we examine 1) solution quality in terms of constraint satisfaction and optimality and 2) computational efficiency over existing convex optimization solver. Simulations are run on an unloaded Macbook Pro with Intel Core i5 8259U CPU @ 2.30GHz.


\textbf{Training Set Generation.} \quad We evaluate our approach on 4 different size of network flow problems \eqref{eqn:network} consisting $20, 50, 200, 1000$ nodes with randomly generated connectivity, and 2 power system benchmarks: the IEEE 14-Bus and 39-Bus system~\cite{athay1979practical}. For generic network flow problems with random topology, it is hard to find specialized algorithms and the best performances are given by commercial solvers~\cite{Mittelmann17}. Similarly, the industry standard for power system problems is to use general purpose solvers~\cite{gomez2018electric}. To generate the training set, we use CVXPY \cite{diamond2016cvxpy} powered by a CVXOPT solver~\cite{andersen2013cvxopt} to solve \eqref{eqn:network} and \eqref{eqn:opf_DC}.  For each network model, we generate $\boldsymbol\ell$ by sampling from uniform distribution. The largest deviation from the nominal are $30\%$ and $80\%$ respectively for network flow and OPF cases respectively. We solve $60,000$ data samples using CVXPY for each network model under each variance setting, and split $20\%$ of the data as test samples.  Overview of evaluated tasks is listed in Table \ref{table:problem} in Appendix \ref{appendix:simulation}.

\textbf{Baselines} \quad We train and compare three other learning models in terms of finding optimal solutions while satisfying all constraints :

\emph{Nearest neighbor for active constraint sets:} This benchmark is a sanity check on whether a deep neural network is needed or a simpler method would suffice. We fit and find a 3-nearest neighbor algorithm achieves the highest classification accuracy in predicting the active constraints.

\emph{End-to-end regression:} Following~\cite{pan2019deepopf}, we construct a $4$-layer neural network to fit the regression task of predicting optimal solution based on input $\boldsymbol\ell$. Mean squared error is used as training loss.

\emph{Classification for active constraint sets: } Following~\cite{ng2018statistical}, we construct a $4$-layer neural network to predict the set of active constraints at optimal solution. We use one-hot encoding for different set of active constraints, and use cross-entropy as the loss function.

For the neural decoder model, the largest model we use is a $4$-layer fully-connected neural network for the 1000-node network flow problem. For other tasks, we find that even a $3$ layer neural network with $200$ neurons on the first layer can produce very competitive results.
We feed load vector $\boldsymbol\ell$  as the input for all above methods. Once the active constraints are predicted, a linear equation solver is used to find final solutions. More simulation setup and architecture details can be found in Appendix \ref{appendix:simulation}.

\begin{table*}[hbt!]
	\centering
	\begin{small}
		\begin{tabular}
			{P{0.2 cm}P{1.9cm} P{0.8cm}P{0.82cm}P{0.8cm}P{0.8cm}P{0.82cm}P{0.8cm}P{0.8cm}P{0.82cm}P{0.8cm}}
			\toprule[0.3mm]
			& &\multicolumn{3}{c}{\textbf{Binding Generators}}&\multicolumn{3}{c}{\textbf{Binding Lines}} &\multicolumn{3}{c}{\textbf{Feasibility}} \\
			\multicolumn{2}{c}{Load Variations}&Low& Meidum& High &Low &Medium&High &Low &Medium&High \\
			\toprule[0.2mm]
			{\multirow{4}{*}{\rotatebox[origin=c]{90}{\textbf{ 39-Bus}}}}&	{Classification} &98.95 &95.28   & 90.19 & 99.24 & 95.96& 86.28 & \textbf{99.99} & 97.38 & 86.51 \\
			&	{End to End} & 93.30& 70.95& 44.38&94.85&88.03&72.22&88.93&54.47&38.75\\
			&	{NearestNeighbor}& 93.40&77.98 &63.90 &59.68 & 78.64 & 71.34 &70.62 & 78.32 &61.29 \\
			&	{NeuralDeco}& \textbf{99.91}&\textbf{96.39} &\textbf{93.95} & \textbf{99.91}& \textbf{96.33} & \textbf{93.48} & \textbf{99.98}& \textbf{99.44}& \textbf{95.75}\\
			\bottomrule[0.3mm]
		\end{tabular}
		\caption{Simulation results for 39-bus OPF problem. With low, medium and high degree of input load variation, we compare three methods in the prediction accuracy of binding generators, binding lines and solution feasibility ratio.}
		\label{table:OPF_Acc}
	\end{small}
\end{table*}

\begin{table}[hbt!]
	\centering
	\begin{small}
		\begin{tabular}
			{P{1.9cm} P{1.5cm}P{1.5cm}P{1.5cm}P{1.5cm}P{1.5cm}P{1.5cm}}
			\toprule[0.3mm]
			\multicolumn{7}{c}{Neural Decoder computation time as a fraction of solvers} \\
			Load variation & NF-20&NF-50&NF-200&NF-1000&OPF-14&OPF-39\\
			\toprule[0.2mm]
			{CVXOPT} & $5.55\%$& $5.30\%$ & $11.58\%$ & $6.59\%$& $7.03\%$ & $10.61\%$  \\
			{CPLEX}& $6.52\%$ & $8.42\%$ & $13.11\%$ & $15.15\%$ & $9.49\%$ & $5.47\%$\\
			{Gurobi}& $6.82\%$ & $7.62\%$ & $16.43\%$ & $18.95\%$& $10.92\%$ & $15.44\%$ \\
			\bottomrule[0.3mm]
		\end{tabular}
		\caption{Average solving time per instance for different solvers.}
		\label{table:Computation_time}
	\end{small}
\end{table}

\textbf{Experimental Results} \quad In Fig.~\ref{fig:NF} we show the probability of encountering infeasible solutions on testing data of network flow problem \eqref{eqn:network}. Infeasible solutions occur either when there is nodal imbalance\footnote{We allow a $5\%$ mismatch when evaluating each node's equality constraint.}, or solved line flows are outside of line capacity bounds. Note that we keep the size of the dataset to be 60,000 to evaluate  generalization properties. For moderate and large networks, these data points are not sufficient to cover the input space and some test samples that do not reside in the same region as the training samples. The proposed Neural Decoder generalizes well to larger networks and has the lowest infeasibility for all test settings. For feasible test instances, the mean solution costs compared to optimal solutions provided by CVXPY is within $0.5\%$.

The other methods perform much poorer than the proposed Neural Decoder. The nearest neighbor approach is not able to learn the correct active constraints, showing that more sophisticated learning approaches are needed. The end-to-end prediction is hard to generalize to larger network and higher variances, mainly because it does not explicitly include line flow constraints, so it tends to produce infeasible flow solutions.
The classification approach is also not a proper learning strategy for larger-scale optimization problems, since the growing number of possible combination of active constraints leads to huge one-hot encodings at the classifier's output. What is more, for unseen combination of active constraints, the classifier is not able to predict the active constraints correctly.

Results on learning for the 39-bus OPF problem are listed in Table~\ref{table:OPF_Acc}, where we report the accuracy of active generators' constraints and lines' constraints separately. For all the feasible solutions provided by Neural Decoder on test samples, the mean cost increase compared to solvers' solution is within $0.7\%$. We notice that proposed method can provide solutions with lower percentage of infeasibility than the sum of error on active generator and line constraints predictions. This is because identification steps of active generators and active lines are performed in sequence and we are making use of the information on total number of active constraints, so the error on each step does not compound and lead to infeasible solutions. 

Compared to the solving process of optimization solvers, proposed neural decoder provides significant speed-up in all testing benchmarks as shown in Table \ref{table:Computation_time}. To make a fair comparison, we use the "core" computation times. That is, using line profiling, we discount the overhead in problem conversion, constraint translation and so on, and only record the time used by lower level linear algebra packages. More results including overhead and statistics are given in Appendix~\ref{app:experiments}. Compared to state-of-the-art solvers, \emph{we can achieve speedup factor of 15 for smaller networks, and for larger networks we can still achieve a 5 fold gain}. We are using standard Python packages for neural network derivations and equation solvers, while further acceleration can be achieved for the neural decoder such as taking a batch of evaluating samples to calculate $\nabla_{\boldsymbol \ell} g_{\theta}(\boldsymbol \ell)$,  or adopting special linear equation solvers to solve the resulting sparse linear equations once all active constraints are identified.

%% file: Conclusion.tex
\vspace{-3pt}
In this work, we presented a learning framework to solve a family of network flow problems by incorporating the structure and parameters of the optimization model. Using neural networks along with their derivatives, we design a decoding strategy to fulfill the task of predicting active constraints and finding constraint-satisfying optimization problem solutions. The method is error-correcting and achieves an order of magnitude speed acceleration compared with current state-of-the-art LP solvers. It also provides significant higher quality solutions in terms of optimality and constraint satisfaction compared to other learning approaches. We hope by incorporating model and parameter information from optimization programs, such fast and robust learning framework can be extended to several future directions, such as learning to solve nonlinear convex optimization problems of specific structures, as well as solving robust and online optimization problems.

%% file: appendix.tex

\section{Modeling Cycle Constraints}
\label{app:cycle}
In power systems, a transmission line is characterized by its admittance (reciprocal of its impedance). The voltage at each of the nodes is described by a complex number. In this paper, we adopt the commonly used DC power flow approximation, which assumes that all voltage magnitudes are fixed, the power flow are driven by the angle differences and the transmission lines are lossless. In particular, if the admittance between nodes $i$ and $k$ is $b_{i,k}$, the power flow $f_{i,k}$ from $i$ to $k$ is $b_{i,k}(\theta_i-\theta_k)$. Consider the nodes $1,\dots, m$ in a cycle of length $m$. Then $\sum_{j=1}^{m} \frac{f_{j,j+1}}{b_{j,j+1}}=\sum_{j=1}^{m} \frac{b_{i,k}(\theta_i-\theta_k)}{b_{j,j+1}}= \sum_{j=1}^{m}(\theta_j-\theta_{j+1})=0$, since all of the angles cancels out. Therefore, a weighted sum of the flows along a cycle must sum up to $0$. For example, in Fig.~\ref{fig:model}, the flows on each link can be changed independently to each other for the network on the left. However, for the power system on the right of Fig.~\ref{fig:model}, (a linear combination of) flows are constrained to sum to $0$.

As described in the main text, due to the existence of  these \emph{cycle constraints}, the edge flows must lie in a subspace in $\mathbb {R}^m$. This space is the so-called fundamental flow space, with a dimension of $n-1$ corresponding to $n-1$ linearly independent equality constraints. The rest of the $m-n+1$ flows are uniquely determined by the cycle constraints. There are many ways to construct this space. A simple construction is to select a spanning tree of the network, where the flows on the edges in this tree are fundamental and form the cycle basis of the network.

\section{Proofs}
\label{appendix:Proof}

\subsection{Proof of Theorem \ref{thm:mu}}
\label{proof:mu}
\begin{proof}
		The proof follows from standard duality theory. We require feasibility of the input $\boldsymbol\ell$ to rule out the dual being unbounded. The Lagrangian of \eqref{eqn:opf_DC} is given by
	\begin{equation}
	\label{eqn:OPF_Lag}
	\begin{split}
	L(\boldsymbol\ell,\bd x, \bd f, \bm \mu,\bar{\bm \lambda}, \underline{\bm \lambda}, \bar{\bm \nu},\underline{\bm \nu}) =\bd c^T \bd x+  \bm \mu^T(\boldsymbol\ell-\bd x-{\tilde{\bd A}} \bd f)+\\
	\bar{\bm \lambda}^T(\bd K \bd f-\bar{\bd f})-\underline{\bm \lambda}^T(\bd K \bd f + \underline{\bd f})-
	\underline{\bm \nu}^T \bd x+\bar{\bm \nu}^T(\bd x-\bar{\bd x})\vspace{-8pt}
	\end{split}
	\end{equation}
	where $\bar{\bm \nu}$ and $\underline{\bm \nu}$ are the dual variables associated with the capacity constraints in \eqref{eqn:opf_DC_gen}, $\bar{\bm \lambda}$ and $\underline{\bm \lambda}$ are the dual variables associated flow capacities in~\eqref{eqn:opf_DC_flow}, and $\bm \mu$ is the dual variable associated with the equality constraint~\eqref{eqn:opf_DC_equ}. The dual variables $\bar{\bm \nu}$, $\underline{\bm \nu}$, $\bar{\bm \lambda}$ and $\underline{\bm \lambda}$ are nonnegative.

	If $\bd x^*, \bd f^*,\bm \mu^*,\bar{\bm \lambda}^*, \underline{\bm \lambda}^*, \bar{\bm \nu}^*,\underline{\bm \nu}^*$ are the optimal primal and dual solutions, then by strong duality
	\begin{equation*}
		J^*(\boldsymbol\ell) = L(\boldsymbol\ell, \bd x^*, \bd f^*,\bm \mu^*,\bar{\bm \lambda}^*, \underline{\bm \lambda}^*, \bar{\bm \nu}^*,\underline{\bm \nu}^*)\vspace{-3pt}
	\end{equation*}
	and differentiating \eqref{eqn:OPF_Lag} gives $\nabla_{\boldsymbol\ell} J^*=\bm \mu^*$.

	We adopt an economics argument first and a rigorous proof through the KKT conditions is given in the next paragraph. Since $\mu_i^*=\frac{\partial J^*}{\partial l_i}$, we can interpret it as the marginal cost of producing one more unit of good at node $i$. If $\mu_i^*<c_i$, the marginal cost of production is lower than the cost at node $i$, which means that node $i$ has zero production ($x_i^*=0$). Conversely, if $\mu_i^*>c_i$, then node $i$ can not produce more to satisfy the demand at a lower cost, then it must be at its upper bound $\bar{x}_i$.

	More rigorously, we can analyze the relationships between binding inequality constraints and optimal dual variables associated with equality constraints. With known $\bm \mu^*$ and separable constraints \eqref{eqn:OPF_dual_2}\eqref{eqn:OPF_dual_3}, the dual problem is decomposable to the optimization problem consisting $\bar{\bm \lambda},\underline{\bm \lambda}$ and $ \bar{ \bm \nu},\underline{\bm \nu}$ respectively. The optimization problem involving $\bar{ \bm \nu},\underline{\bm \nu}$ can be reformulated as
	\begin{subequations}
		\label{eqn:OPF_dual_nu}
		\begin{align}
		\min_{\bar{ \bm \nu},\underline{\bm \nu}} \quad  & \bar{\bm \nu}^T\bar{\bd x} \label{eqn:OPF_nu_1}\\
		\text { s.t. } \quad &\bar{\bm \nu}-\underline{\bm \nu}=\bm \mu^*-\bd c,   \label{eqn:OPF_nu_2}\\
		& \bar{ \bm \nu}\geq \bd 0,\, \underline{\bm \nu}\geq \bd 0
		\end{align}
	\end{subequations}
	which can read out the result without solving explicitly. Given that $\bar{ \bm \nu},\underline{\bm \nu}\geq \bd 0$, for node $i$, if $\mu_i^*-c_i>0$, $\bar{ \nu}_i>0,\underline{\nu}_i=0$, and $x_i=\bar{x}_i$ ; if $\mu_i^*-c_i<0$, $\bar{ \nu}_i=0,\underline{\nu}_i>0$, and  $x_i=0$;  if $\mu_i^*-c_i=0$, $\bar{ \nu}_i=\underline{\nu}_i=0$ and corresponding $x_i$ is not binding.
\end{proof}

\subsection{Proof of Lemma~\ref{lem:upsilon}} \label{app:upsilon}
\begin{proof}
The original optimization problem is reproduced below:
\begin{subequations}
	\begin{align*}
	\max_{\bar{{\bm \lambda}},\underline{\bm \lambda}} \quad  &-\underline{\bm \lambda}^T\underline{\bd f}-\bar{\bm \lambda}^T\bar{\bd f}\\
	\text { s.t. } \;  &  \bd K^T(\bar{\bm \lambda}-\underline{\bm \lambda})=\tilde{\bd A}^T\bm \mu^* \\
	& \bar{\bm \lambda}\geq \bm 0,   \,    \underline{\bm \lambda}\geq  \bm 0.
	\end{align*}
\end{subequations}
Assuming $\underline{\bd f}=\bar{\bd f}$ and let $\bm \upsilon= \bar{\bd f} \odot \bar{\bm \lambda}+\underline{\bd f} \odot \underline{\bm \lambda}$ and $\bd v=\bar{\bm \lambda}-\underline{\bm \lambda}$, where $\odot$ is componentwise multiplication. Then the optimization problem becomes
\begin{subequations}
	\begin{align*}
	\min_{\bm \upsilon, \bd v} \quad & \sum_{i=1}^n \upsilon_i \\
	\text { s.t. } \;  &  \bd K^T \bd v=\tilde{\bd A}^T\bm \mu^* \\
	& \bm \upsilon \geq \bar{\bd f} \odot \bd v \\
	& \bm \upsilon \geq -\bar{\bd f} \odot \bd v,
\end{align*}
\end{subequations}
where the last two inequalities come from the nonnegativity constraint of $\underline{\bm \lambda}$. Suppose $\bm \upsilon, \bd v$ are optimal solutions. Because we are minimizing the sum of the components of $\bm \upsilon$, $\upsilon_i=\bar{f}_i \max(v_i,-v_i)$. This is equivalent to $\upsilon_i=\bar{f}_i |v_i|$.
\end{proof}

\subsection{Proof of Lemma \ref{lem:error}}
\label{app:error}
We use the following lemma:
\begin{lemma}
	\label{lemma:region}
For an LP problem $\{ \bd x^*=\arg \min_{\bd x} \bd c^T \bd x | \bd {Ax}= \bd b, \bd x \geq  \bd 0 \}$ with $\bm A \in \mathbb{R}^{m\times n}, \, \bd x\in \mathbb{R}^m, m<n$, let $(1),(2),...,(m)$ be the indices of selected columns of $\bd A$, such that $\bm B \bm x^*= \bd b$ with $\bd B=[ \bd a_{(1)} \, ... \, \bd a_{(m)}  ]\in \mathbb{R}^{(m\times m)}$ as an invertible basis. 	For any $\tilde{\bd b}=\bd b + \bm \delta$, the optimal solution is still given by $\bd B$ if and only if $\bm B^{-1}\bd b+ \bm B^{-1}\bm \delta \geq  \bd 0$.
\end{lemma}
The proof of Lemma~\ref{lem:error} follows from Lemma~\ref{lemma:region} by converting our LP of interest to the standard form. Despite the latter being a known result in linear programming, we have not been able to find the proof of the vector form in existing literature (the component by component result can be found in~\cite{bertsimas1997introduction, jansen1997sensitivity}). Therefore we provide the following proof for completeness.
\begin{proof}
	To find the region where $\tilde{\bd b}=\bd b+\bm \delta$ has the same set of active constraints at the optimal solution, denote the new optimal solution as $\tilde{\bd x}^*$ that satisfies $\bd {A\tilde{ x}}=\tilde{\bd b}$. So the new optimization problem involving $\tilde{\bd b}$ becomes
	\begin{subequations}
		\begin{align}
		\tilde{\bd x}=\arg \min \quad &\bd c^T \bd x\\
		s.t. \quad & \bd A \bd x=\bd b+ \bm \delta\\
		& \bd x\geq \bd 0
		\end{align}
	\end{subequations}

	Since we have $\bd x^*=\bd B^{-1}\bd b \geq 0$ with input $\bd b$, and since $\bd x^*$ and $\tilde{\bd x}$ can be represented by the same basis $\bd B \in \mathbb{R}^{m\times m}$, we have
	\begin{equation}
	\label{eqn:appen_1}
	\bd B^{-1}(\bd b + \bm \delta)\geq \bd 0
	\end{equation}
 to ensure the optimal solution's feasibility. So the resulting $\bm \delta$ must satisfy \eqref{eqn:appen_1}.

 On the other hand, if $\bd B^{-1}(\bd b + \bm \delta)\geq 0$, while $\tilde{ \bd x}=\bd B^{-1}(\bd b + \bm \delta)$ satisfies both equality and inequality constraints. By checking the KKT conditions, it is also the optimal solution, which completes the proof.
\end{proof}

\section{Dictionary Learning}
\label{app:dictionary}

Lemma \ref{lem:error} also provides a fast alternative for identifying $\bar{ \bm \lambda}, \underline{\bm \lambda}$ and associated active constraints in \eqref{eqn:OPF_line_dual}~\cite{deng2016probabilistic}. We now discuss how to make use of topology and known parameters of the optimization problem to accelerate the robust decoding process of active constraints identification.  During training or data generation process, we can keep an offline, finite-length dictionary of unique $\mathpzc{M}=\{ \tilde{\bd A}^T \bm \mu^{*}_{(k)}\}$ along with the set of optimal dual variables $\{[\bar{\bm{\lambda}}^{*}_{(k)}, \underline{\bm{\lambda}}^{*}_{(k)}]\}$, $\mathpzc{P}_{\tilde{\bd A}^T  \bm \mu^*_{(k)}}$ and active constraints,  respectively. When $\boldsymbol\ell$ comes along with neural network prediction $g_\theta(\boldsymbol\ell)$ at testing time, with $\bm \delta_{(k)} = \tilde{\bd A} \nabla_{\boldsymbol\ell} g_\theta(\boldsymbol\ell) - \tilde{\bd A}^T \bm \mu^{*}_{(k)}$, one lookup operation is required to determine if $\bm \delta_{(k)} \in \mathpzc{P}_{\tilde{\bd A}^T \bm \mu^*_{(k)}}$, and we can directly identify the value of $\bm \lambda$. If the dictionary does not include the region of active constraint sets, $\mathpzc{L}_1$ minimization problem \eqref{eqn:OPF_line_L1} is solved via IHT algorithm. Such design of active constraints prediction contains an error-correcting step by design, where a relatively inaccurate $\nabla_{\boldsymbol\ell} g_\theta(\boldsymbol\ell)$ (compared to $\bm \mu^*$) can still lead to the correct set of binding constraints.

\section{3-bus Example}
\begin{figure}[t]
 	\begin{minipage}[c]{0.55\textwidth}
 		\centering
 		\includegraphics[width=1.0\linewidth]{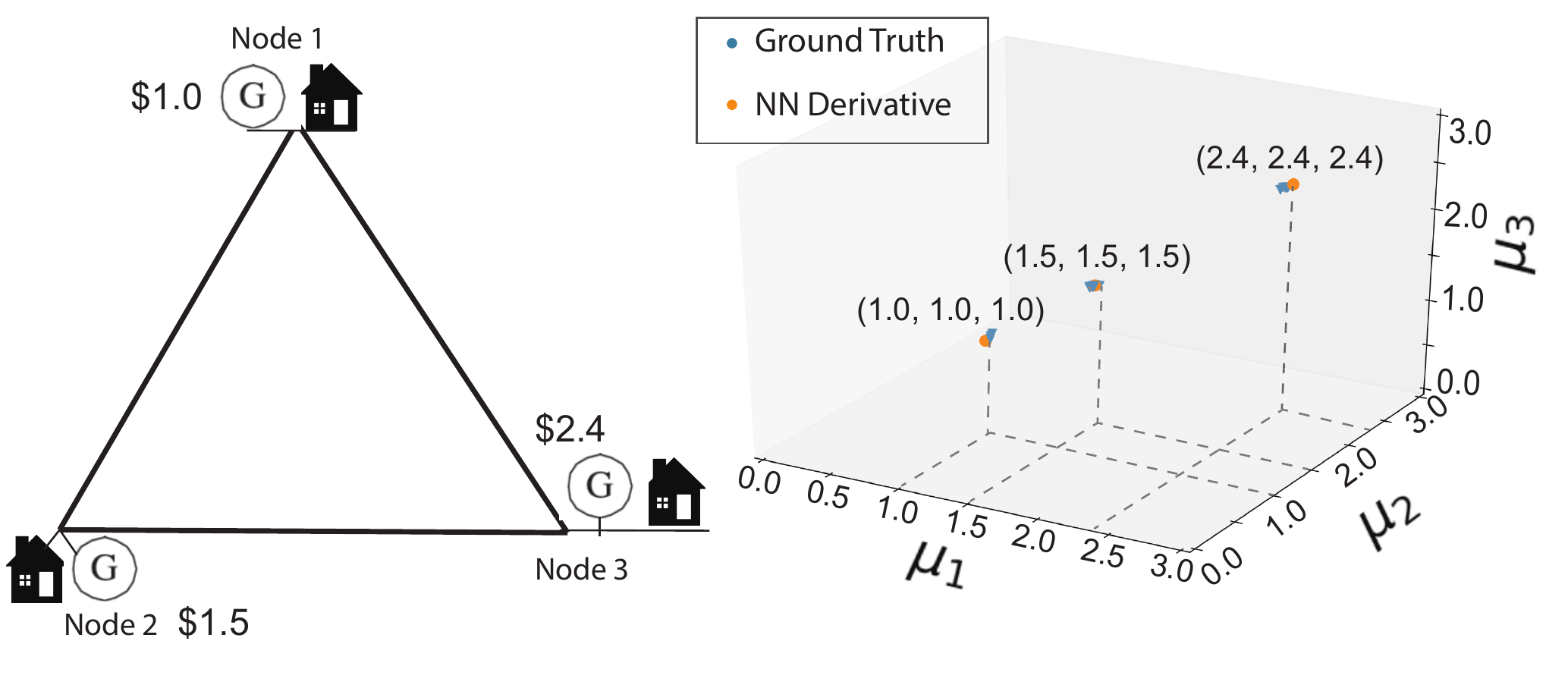}
 	\end{minipage}\hfill
 	\begin{minipage}[c]{0.4\textwidth}
 		\caption{OPF problem for a 3-bus system along with the visualization for dual variables $\mu$. (Left) each bus connect with a dispatchable generator and a variable load, (right) and proposed approach learns to predict $\mu$ correctly.}
 		\label{fig:3bus}
 	\end{minipage}
 \end{figure}

 In Fig.\ref{fig:3bus}, we also show a toy example of proposed learning approach, which includes a 3-bus OPF example with $c_i=1.0, 1.5, 2.4$ respectively. There are 5 optimization variables (3 generators, 2 line flows), 6 inequality constraints associated with generators and 6 inequality constraints associated with line flows. We generate $100$ test load samples which are generated from a uniform distribution, and visualize $\bm \mu^*$. Our learning approach accurately predicts $\bm \mu^*$, and is then able to identify binding generation constraints.
 Even if learned $\bm \mu^*$ has errors, our approach is still robust and can learn the correct set of active constraints associated with the optimal solution.

\section{Main Algorithm} \label{app:algorithm}
In Algorithm \ref{alg:algorithm}, we summarize our learning procedures described in the main texts along with the dictionary learning discussed in Appendix \ref{app:dictionary}. For each testing instance $\boldsymbol{\ell}$, we firstly calculate $\nabla_{\boldsymbol \ell} g_\theta(\boldsymbol{\ell})$, and find the binding nodal constraints based on Theorem \ref{thm:mu}. If there is not an offline-collected dictionary or the resulting $\bm \delta_{(j)}$ is not included in the offine dictionary, then an $\mathpzc{L}_1$ minimization step is implemented to identify active line flow constraints. Otherwise we can directly find the set of constraints by looking up the codeword in the dictionary. A final linear equation solver is utilized to find the optimal solutions using identified active inequality constraints along with the equality constraints of the optimization problem. The major computation steps include neural network gradient calculation, the $\mathpzc{L}_1$ or dictionary lookup, and the linear equaltions solving.

\renewcommand{\algorithmicrequire}{\textbf{Input:}}
\renewcommand{\algorithmicensure}{\textbf{Parameters:}}
\newcommand{\algorithmicbreak}{\textbf{break}}
\newcommand{\BREAK}{\STATE \algorithmicbreak}

\begin{figure}[t]
	\begin{algorithm}[H]
		\caption{ Neural Decoder for Active Constraints}
		\label{alg:algorithm}
		\begin{algorithmic}[1]
			{\small
				\REQUIRE $\boldsymbol\ell$, trained model $g_{\theta}$, number of active constraints $K$
				\REQUIRE (\textbf{Optional}, see Appendix \ref{app:dictionary}) Dictionary  $\{ \tilde{\bd A}^T \bm \mu^{*}_{(k)}\}$, $\mathpzc{P}_{\tilde{\bd A}^T \bm \mu^*_{(k)}}\, ,\bd f_{(k)}$
				\ENSURE $\epsilon >0$, Optimization model $\tilde{\bd A}, \bar{\bd f}, \underline{\bd f},  \bd c$, $FLAG=0$
				\STATE Find $\nabla_{\boldsymbol \ell} g_\theta(\boldsymbol{\ell})$

				\FOR{$i=1,...,n$}
				\IF{ $\nabla_{l_i} g_\theta(\boldsymbol{\ell})-c_i<-\epsilon$}
				\STATE {$x_i=0$, $K=K-1$}
				\ELSIF {$\nabla_{l_i} g_\theta(\boldsymbol{\ell})-c_i > \epsilon$}
				\STATE $x_i=\bar{x}_i$, $K=K-1$
				\ENDIF
				\ENDFOR
				\FOR{$j=1,...,|M|$}
				\STATE {$\bm \delta_{(j)} = \tilde{\bd A} \nabla_{\boldsymbol\ell} g_\theta(\boldsymbol\ell) - \tilde{\bd A}^T \bm \mu^{*}_{(j)}$}
				\IF {$\bm \delta_{(j)} \in\mathpzc{P}_{\tilde{\bd A}^T \bm \mu^*_{(j)}}$}
				\STATE {Identify active flow constraints $\bd f_{(j)}$, $FLAG=1$}
				\BREAK
				\ENDIF
				\ENDFOR
				\IF {$FLAG=0$}
				\STATE {\texttt{IHTSolve}(\eqref{eqn:OPF_line_L1}, sparsity=$K$) }
				\ENDIF
				\STATE {\texttt{EquationSolve}($\bd x+\bd A\bd f=\boldsymbol{\ell}$, active constraints set)}
			}
		\end{algorithmic}
	\end{algorithm}
\end{figure}

\begin{table}
	\centering
	\begin{small}
		\begin{tabular}
			{m{2.3cm}  P{2.3cm}P{2.3cm}P{2.6cm}}
			\toprule[0.3mm]
			Problem Instance &Nodal Constraints &Line Constraints &Problem Input\\
			\toprule[0.2mm]
			Network Flow &up to 2000& up to 3394& nodal net demands  \\
			OPF IEEE 14-Bus& up to 24& up to 47 & nodal loads\\
			OPF IEEE 39-Bus &up to 59 & up to 100& nodal loads \\
			
			\bottomrule[0.3mm]
			
		\end{tabular}
		\caption{Experiments Overview. }
		\label{table:problem}
	\end{small}
\end{table}

\section{Details of Experiments} \label{app:experiments}

\subsection{Modeling of Network Flow Problem}
We note that formulation showed in \eqref{eqn:network} is a general formulation of \emph{minimum-cost network flow problem}. For instance, in the transportation problem, we could treat the nodes either sources or sinks of the product, and thus we can interpret the nodal equality as $\sum_{k=1}^{n}f_{jk}=a_j$ for all source nodes $j (j=1,2,...,m)$ and $\sum_{j=1}^{m}f_{jk}=d_k$ for all destination nodes $k\; (k=1,2,...,n)$. In the shortest-path problem, $C_i$ denotes the distance associated with each edge $f_i$. With a net supply of one unit at the source, we have $\sum_{k=1}^{n}f_{jk}=1$  for the source node $j$, and $\sum_{j=1}^{m}f_{jk}=1$  for the destination node $k$, and the objective is to send one unit of flow from the source to the link at minimum cost.

\begin{table*}[hbt!]
	\centering
	\begin{small}
		\begin{tabular}
			{P{0.2 cm}P{1.9cm} P{0.8cm}P{0.82cm}P{0.8cm}P{0.8cm}P{0.82cm}P{0.8cm}P{0.8cm}P{0.82cm}P{0.8cm}}
			\toprule[0.3mm]
			& &\multicolumn{3}{c}{\textbf{Binding Generators}}&\multicolumn{3}{c}{\textbf{Binding Lines}} &\multicolumn{3}{c}{\textbf{Feasibility}} \\
			\multicolumn{2}{c}{Load Variations}&Low& Meidum& High &Low &Medium&High &Low &Medium&High \\
			\toprule[0.2mm]
			{\multirow{4}{*}{\rotatebox[origin=c]{90}{\textbf{ 14-Bus}}}}& {Classification} &94.64 &78.07   &71.37 & 98.52 & 98.06& 94.07& 93.37& 79.08&75.38 \\
			&	{End to End} &96.36 & 96.04   &90.39& 96.06 & 95.78& 93.77& 96.64& 96.14 &86.31 \\
			&	{NearestNeighbor}& 89.03&88.92&73.65&92.44&89.68&84.33&88.72&81.83&70.81\\
			&	{NeuralDeco}& \textbf{98.49}&\textbf{96.15} &\textbf{93.19} & \textbf{99.91}& \textbf{98.16} & \textbf{93.47} & \textbf{99.23}& \textbf{98.14}& \textbf{95.97}\\
			\bottomrule[0.3mm]
		\end{tabular}
		\caption{Simulation results for OPF 14-bus problem. With varying input loads (low, medium and high), we compare four methods in the prediction performance of  binding generators (accuracy), binding lines (accuracy) and solution feasibility ratio. }
		\label{table:OPF_Acc_14}
	\end{small}
\end{table*}

\begin{figure}
	\centering
	\includegraphics[width=0.99\linewidth]{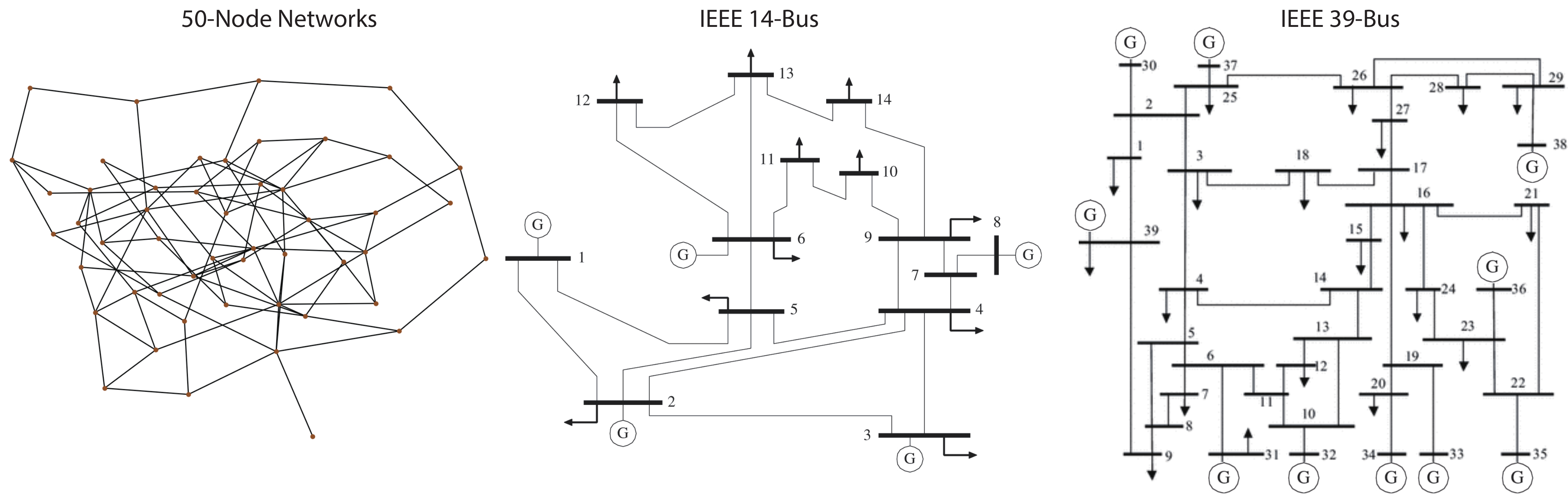}
	\caption{\footnotesize The example network for our simulations on network flow and OPF problem.}
	\label{fig:graph}
\end{figure}

\subsection{Simulation Details and Ablation Study}
\label{appendix:simulation}
In the data generation process, we collected the input data, the output optimal costs, the dual variables and optimal solutions for each network. We collected all the data samples which are feasible for optimization problem \eqref{eqn:network} and \eqref{eqn:opf_DC}. With a nominal load level $\bar{\boldsymbol\ell}$, in the standard network flow problem, we generate input samples from the uniform distribution $[0.9\boldsymbol\ell,\, 1.1\boldsymbol\ell]$, $[0.8\boldsymbol\ell,\, 1.2\boldsymbol\ell]$, $[0.7\boldsymbol\ell,\, 1.3\boldsymbol\ell]$; in the OPF problem, we generate input samples from the uniform distribution $[0.8\boldsymbol\ell,\, 1.2\boldsymbol\ell]$, $[0.5\boldsymbol\ell,\, 1.5\boldsymbol\ell]$, $[0.2\boldsymbol\ell,\, 1.8\boldsymbol\ell]$.

\begin{figure}[t]
	\centering
	\includegraphics[width=0.79\linewidth]{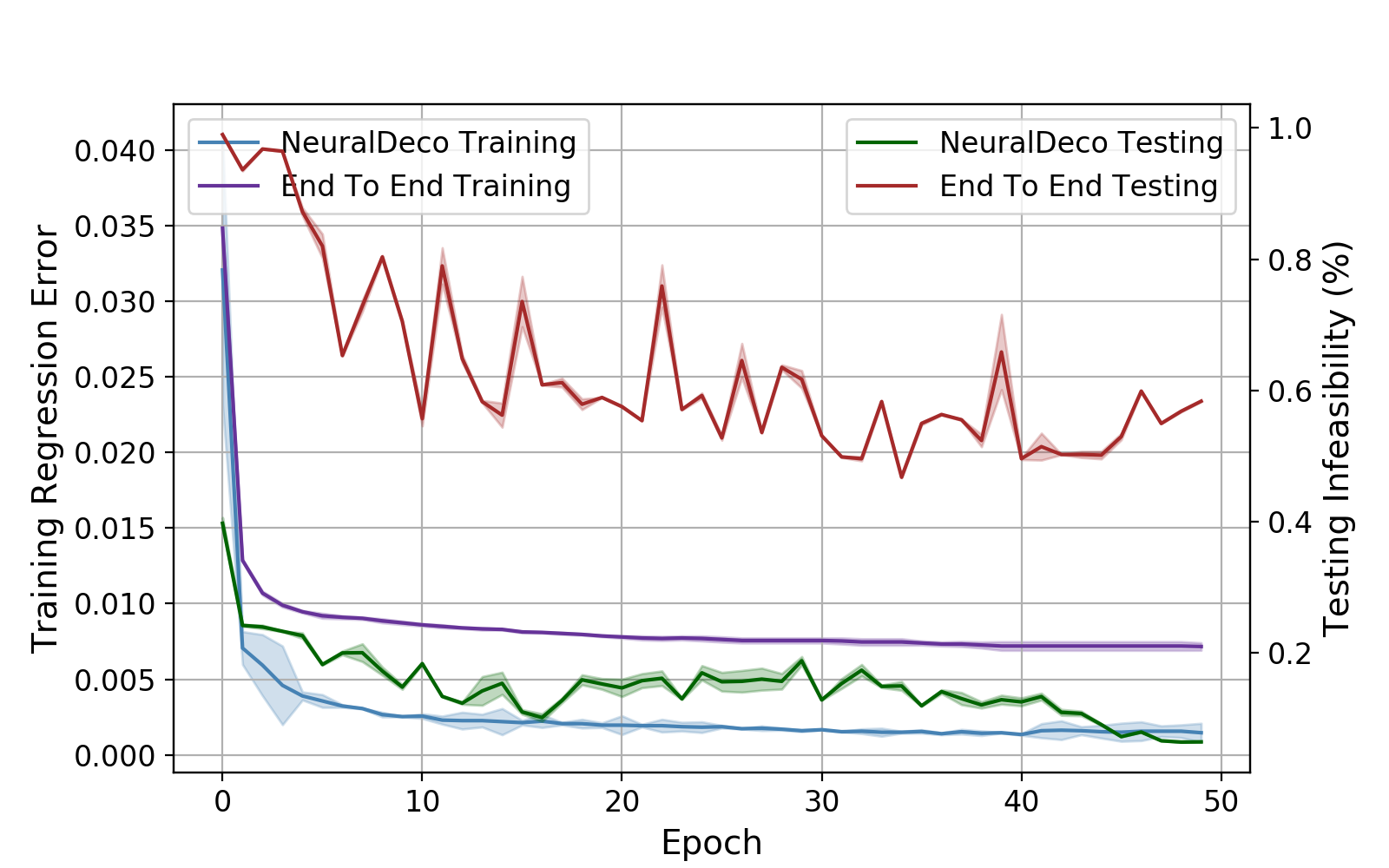}
	\caption{\footnotesize Training error and testing solution feasibility for IEEE 39-bus case. Mean and standard deviation are shown for $3$ runs. }
	\label{fig:OPF}
\end{figure}

\begin{figure}[ht]
	\centering
	\includegraphics[width=0.7\linewidth]{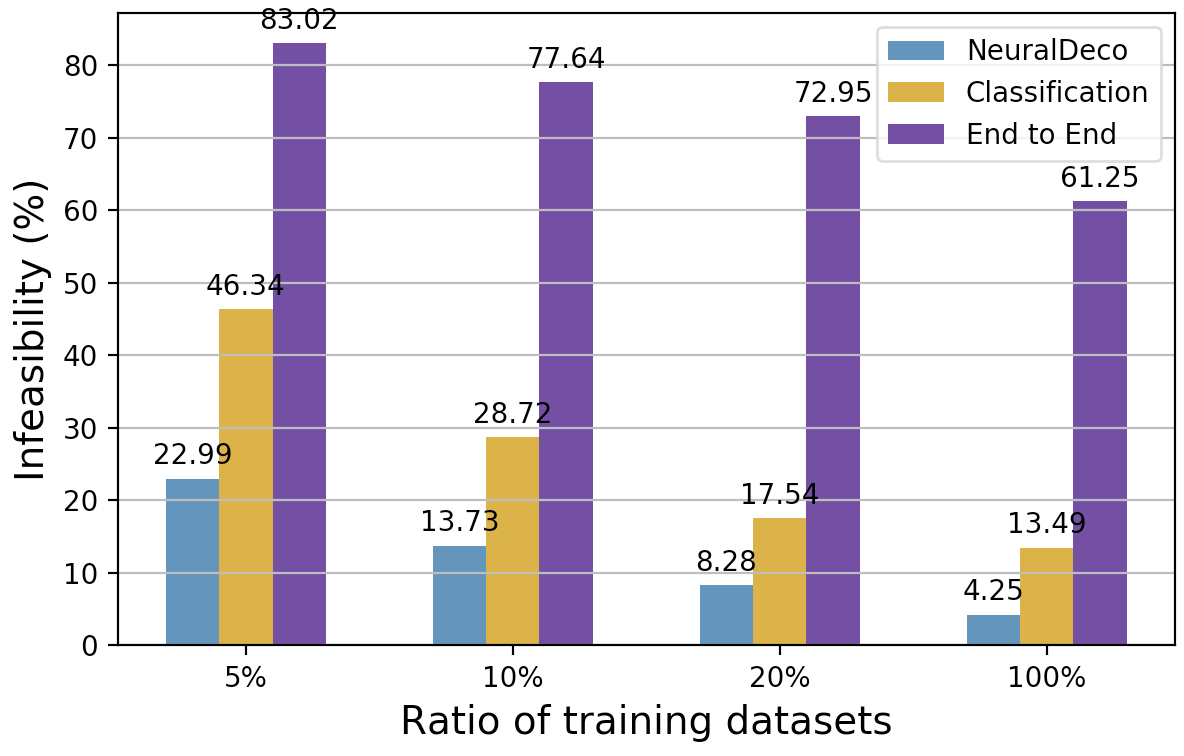}
	\caption{\footnotesize The solution infeasibility of testing instances with varying size of training datasets on the OPF 39-bus case with high input variations. }
	\label{fig:OPF_data}
\end{figure}

In Figure \ref{fig:graph}, we show the underlying graph for network flow and OPF simulation. We generated random connected graph for the network flow problem. The IEEE 14-bus system is known as an approximation of America Electric Power system, while 39-bus system is also known as 10-machine New England Power System, both representing the practical operations of the electricity grid.

For the network flow 1000-bus problem, we use a 4-layer feedforward neural network with $\{1000, 500, 500, 300\}$ neurons at each layer with ReLU activation function. For the OPF IEEE 39-bus problem, we use a 3-layer network with $\{100, 30, 20\}$ neurons at each layer. Based on our evaluations, we found the gains to be negligible when using more number of layers or nodes under the current simulation setting. For the practical implementation of Hamming distance loss between predicted active constraints and actual active constraints,  we only penalize the wrong predictions using value of $\nabla_{\boldsymbol\ell} g(\theta) \boldsymbol\ell - \bd c$ compared to each constraint's binding status.

When reporting the results on testing instances, we treat the following solutions as ``infeasible":
\begin{itemize}
	\item \emph{Infeasible Nodal Values:} When the nodal inequality does not hold for any of the node;
	\item \emph{Infeasible Sources:} When the solution for any of  sources' output or generator's generations are out of the capacity limits;
	\item \emph{Infeasible Line Flows:}  When the solution for any of lines' flow are out of line capacity.
\end{itemize}

We allow $5\%$ nodal mismatch for equality constraints or $5\%$ out of bounds for inequality constraints when evaluating the solutions of different machine learning models. Otherwise solutions are treated as infeasible to the optimization problems.

Results on learning for OPF problem for  IEEE 14-bus system are listed in Table~\ref{table:OPF_Acc_14}. Nearest neighbor approach finds the nearest neighbors on the Euclidean space, which is a mismatch from the decision boundary of the set of active constraints. Thus even on the test cases with low input variations and high similarity between input instances, nearest neighbor search can identify wrong set of binding constraints, leading to high solution infeasibility ratio.

We also investiage why learning the optimal costs rather than end-to-end learning the optimal solutions are a more appealing approach. In Figure \ref{fig:OPF}, we visualize the training progress during 50 training epochs on the OPF task for IEEE 39-bus system. We also evaluate model predictions' feasibility on test dataset after each epoch. Both the regression loss for end to end training and Neural Decoder training converges quickly to loss values smaller than $0.01$. However, after the few epochs' progress on increasing the solution feasibility on testing data, the end-to-end predictions fail to provide practical, feasible solutions. Solution  feasibility of end-to-end approach also varies significantly for different epochs, indicating that the mean squared error defined on solution prediction is not a good indicator for solution feasibility.

In Figure \ref{fig:OPF_data}, we investigate the performance of learning algorithms with smaller size of training datasets. With only $5\%$ of original training data, proposed neural decoder achieves much better accuracy on active constraints identification as well as feasible solutions prediction. Due to the fact that our approach is error-correcting even when training data is limited and there is prediction error, our decoding process can generalize better with small training dataset. This may further help the practical decision making process when solution data collection is limited.

\begin{table*}[!t]
	\centering
	\begin{small}
		\begin{tabular}
			{P{1.7cm} P{1.11cm}P{1.11cm}P{1.01cm}P{1.01cm}P{1.01cm}P{1.11cm}P{1.11cm}P{1.11cm}P{1.01cm}P{1.01cm}}
			\toprule[0.35mm]
			Time(s)  &\multicolumn{4}{c}{NF 20-Node}&\multicolumn{4}{c}{NF 50-Node}\\

			& NeualDeco& CVXOPT&CPLEX&Gurobi& NeuralDeco&CVXOPT& CPLEX&Gurobi\\
			\toprule[0.2mm]
			{Low} & 0.0008& 0.0152& 0.0114&0.0121&0.0009& 0.0292&0.0166&0.0157 \\
			{Medium}& 0.0006&0.0147 &0.0133 &0.0134 &0.0010& 0.0267&0.0182&0.0153\\
			{High}& 0.0006&0.0149& 0.0139 &0.0125 &0.0010&0.0256&0.0207&0.0168\\
			{Solving Time} &0.0006&0.0108&0.0092&0.0088&0.0008&0.0151&0.0095&0.0105\\
			\toprule[0.2mm]
			&\multicolumn{4}{c}{NF 200-Node}&\multicolumn{4}{c}{NF 1000-Node}\\
			& NeualDeco& CVXOPT&CPLEX&Gurobi& NeuralDeco&CVXOPT& CPLEX&Gurobi\\
			\toprule[0.2mm]
			{Low} & 0.0108 & 0.1698 &0.0805&0.0716&0.1709&2.6837&0.6351&0.5152 \\
			{Medium}& 0.0121 & 0.1622&0.0928&0.0816&0.1727&2.7754&0.7082&0.5279\\
			{High}& 0.0108&0.1698&0.1107 & 0.1062&0.1733& 3.1096 &0.6875& 0.5574\\
			{Solving Time}  & 0.0094&0.0812 & 0.0717&0.0572 & 0.0922 &1.3995 &0.6085& 0.4864\\
			\bottomrule[0.35mm]
		\end{tabular}
		\caption{Computation time per instance for NF problem using CVXOPT, CPLEX, Gurobi and proposed Neural Decoder. }
		\label{table:Computation_time2}
	\end{small}
\end{table*}

\begin{table*}[t]
	\centering
	\begin{small}
		\begin{tabular}
			{P{1.7cm} P{1.11cm}P{1.11cm}P{1.01cm}P{1.01cm}P{1.01cm}P{1.11cm}P{1.11cm}P{1.11cm}P{1.01cm}P{1.01cm}}
			\toprule[0.3mm]
			Time(s)  &\multicolumn{4}{c}{OPF 14-Bus}&\multicolumn{4}{c}{OPF 39-Bus} \\

			&  NeuralDeco&CVXOPT&CPLEX &Gurobi &NeuralDeco &CVXOPT &CPLEX &Gurobi\\
			\toprule[0.2mm]
			{Low} & 0.0013 & 0.0333  & 0.0247&0.0149 &0.0022 & 0.0283 &0.0497&0.0194 \\
			{Medium}&  0.0014&0.0233 & 0.0211& 0.0156&0.0024&0.0285 &0.0531&0.0179\\
			{High}& 0.0015&0.0237& 0.0248&0.0162 &0.0028&0.0291 & 0.0529 &0.0198\\
			{Solving Time}& 0.0013&0.0185& 0.0137& 0.0119&0.0021&0.0198 & 0.0384&0.0136\\
			\bottomrule[0.3mm]

		\end{tabular}
		\caption{Computation time per instance for OPF problem using CVXOPT, Gurobi and CPLEX and proposed Neural Decoder. }
		\label{table:Computation_time3}
	\end{small}
\end{table*}

\begin{figure}[h]
	\centering
	\includegraphics[width=0.75\linewidth]{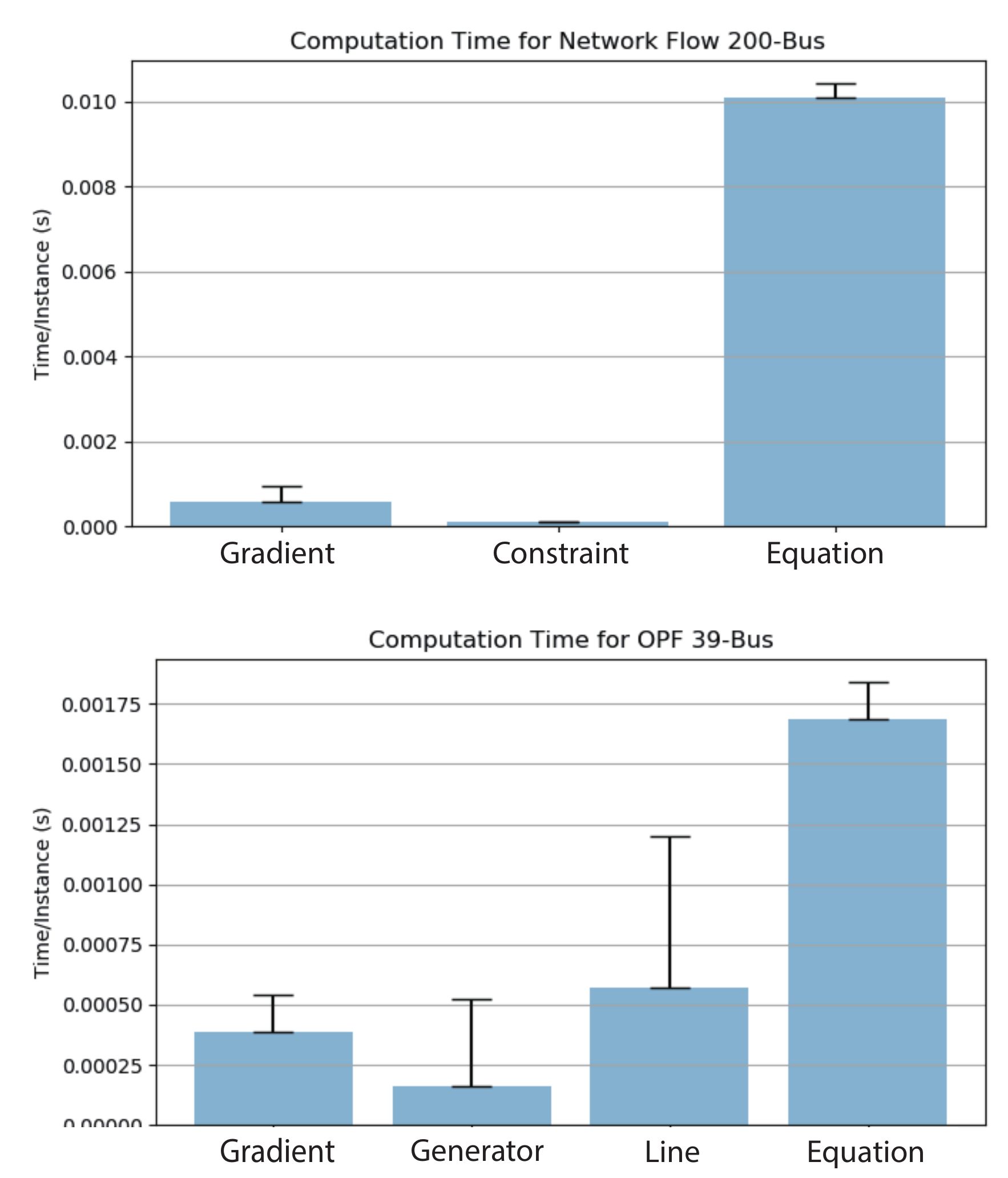}
	\caption{\footnotesize The computation time analysis on network flow problem and optimal power flow problem.}
	\label{fig:time}
\end{figure}

\newpage
\subsection{Computation Time Analysis} As is reported in Table \ref{table:Computation_time2}, proposed Neural Decoder approach could achieve faster computation of optimal solutions for various settings of network flow problems compared to CVXOPT, CPLEX and Gurobi solvers. We note that we are making use of off-the-shelf packages for neural networks differentiation and equation solving. More acceleration may be expected if the specific neural network architecture or the equation structure is considered. Proposed method can also achieve acceleration when batch operattions are adopted to calculate the gradients $\nabla_{\boldsymbol\ell} J^*$, which is a practical implementation procedure when multiple input instances are considered.

Here we report two metrics on computation time. The first is the total runtime of each of the solver. Since we only take a straightforward gradient prediction step when solving the problem, and that removes the need for iterations in optimization solver. The main computation hurdle of our method is actually the linear equation solving step using NumPy packages (Fig. \ref{fig:time}) due to the modeling and interfacing steps with the lower level LAPACK functions. By doing a further inspection on line-profiling analysis, our second metric tries to exclude the overhead and to only include the runtime of the “core” numerical algorithm executions, and the results comparing ``solving time'' are summarized in Table \ref{table:Computation_time2} and Table \ref{table:Computation_time3}, which is also used to calculate the time ratio shown in the main text.

%

%% file: main.bbl
\begin{thebibliography}{10}
\providecommand{\url}[1]{#1}
\csname url@samestyle\endcsname
\providecommand{\newblock}{\relax}
\providecommand{\bibinfo}[2]{#2}
\providecommand{\BIBentrySTDinterwordspacing}{\spaceskip=0pt\relax}
\providecommand{\BIBentryALTinterwordstretchfactor}{4}
\providecommand{\BIBentryALTinterwordspacing}{\spaceskip=\fontdimen2\font plus
\BIBentryALTinterwordstretchfactor\fontdimen3\font minus
  \fontdimen4\font\relax}
\providecommand{\BIBforeignlanguage}[2]{{%
\expandafter\ifx\csname l@#1\endcsname\relax
\typeout{** WARNING: IEEEtran.bst: No hyphenation pattern has been}%
\typeout{** loaded for the language `#1'. Using the pattern for}%
\typeout{** the default language instead.}%
\else
\language=\csname l@#1\endcsname
\fi
#2}}
\providecommand{\BIBdecl}{\relax}
\BIBdecl

\bibitem{bertsimas1997introduction}
D.~Bertsimas and J.~N. Tsitsiklis, \emph{Introduction to linear
  optimization}.\hskip 1em plus 0.5em minus 0.4em\relax Athena Scientific
  Belmont, MA, 1997, vol.~6.

\bibitem{wen2015optimal}
Z.~Wen, D.~O’Neill, and H.~Maei, ``Optimal demand response using device-based
  reinforcement learning,'' \emph{IEEE Transactions on Smart Grid}, vol.~6,
  no.~5, pp. 2312--2324, 2015.

\bibitem{dalal2016hierarchical}
G.~Dalal, E.~Gilboa, and S.~Mannor, ``Hierarchical decision making in
  electricity grid management,'' in \emph{International Conference on Machine
  Learning}, 2016, pp. 2197--2206.

\bibitem{Jimenez19}
J.~Jimenez, ``How the {Texas} power grid braces against rolling blackouts as
  summer heat looms,'' \emph{The Dallas Morning News}, 2019.

\bibitem{Penn19}
I.~Penn, ``This did not go well: Inside {PG}\&{E}’s blackout control room,''
  \emph{The New York Times}, 2019.

\bibitem{glover2012power}
J.~D. Glover, M.~S. Sarma, and T.~Overbye, \emph{Power system analysis \&
  design, SI version}.\hskip 1em plus 0.5em minus 0.4em\relax Cengage Learning,
  2012.

\bibitem{stott2009dc}
B.~Stott, J.~Jardim, and O.~Alsa{\c{c}}, ``Dc power flow revisited,''
  \emph{IEEE Transactions on Power Systems}, vol.~24, no.~3, pp. 1290--1300,
  2009.

\bibitem{deng2016probabilistic}
W.~Deng, Y.~Ji, and L.~Tong, ``Probabilistic forecasting and simulation of
  electricity markets via online dictionary learning,'' \emph{arXiv preprint
  arXiv:1606.07855}, 2016.

\bibitem{pan2019deepopf}
X.~Pan, T.~Zhao, and M.~Chen, ``Deepopf: A deep neural network approach for
  security-constrained dc optimal power flow,'' \emph{arXiv preprint
  arXiv:1910.14448}, 2019.

\bibitem{amos2017optnet}
B.~Amos and J.~Z. Kolter, ``Optnet: Differentiable optimization as a layer in
  neural networks,'' in \emph{Proceedings of the 34th International Conference
  on Machine Learning-Volume 70}.\hskip 1em plus 0.5em minus 0.4em\relax JMLR.
  org, 2017, pp. 136--145.

\bibitem{agrawal2019differentiable}
A.~Agrawal, B.~Amos, S.~Barratt, S.~Boyd, S.~Diamond, and J.~Z. Kolter,
  ``Differentiable convex optimization layers,'' in \emph{Advances in Neural
  Information Processing Systems}, 2019, pp. 9558--9570.

\bibitem{bertsimas2012adaptive}
D.~Bertsimas, E.~Litvinov, X.~A. Sun, J.~Zhao, and T.~Zheng, ``Adaptive robust
  optimization for the security constrained unit commitment problem,''
  \emph{IEEE transactions on power systems}, vol.~28, no.~1, pp. 52--63, 2012.

\bibitem{bixby2015computational}
R.~E. Bixby, ``Computational progress in linear and mixed integer
  programming,'' \emph{Presentation at ICIAM}, vol. 2015, 2015.

\bibitem{Heidel16}
T.~Heidel, F.~Pan, S.~Elbert, C.~DeMarco, and H.~Mittelmann, ``Optimal power
  flow competition design considerations,'' in \emph{Increasing Market and
  Planning Efficiency through Improved Software, FERC Technical Conference},
  2016.

\bibitem{gomez2018electric}
A.~Gomez-Exposito, A.~J. Conejo, and C.~Canizares, \emph{Electric energy
  systems: analysis and operation}.\hskip 1em plus 0.5em minus 0.4em\relax CRC
  press, 2018.

\bibitem{dantzig2003max}
G.~Dantzig and D.~R. Fulkerson, ``On the max flow min cut theorem of
  networks,'' \emph{Linear inequalities and related systems}, vol.~38, pp.
  225--231, 2003.

\bibitem{vanderbei2015linear}
R.~J. Vanderbei \emph{et~al.}, \emph{Linear programming}.\hskip 1em plus 0.5em
  minus 0.4em\relax Springer, 2015.

\bibitem{rao2019engineering}
S.~S. Rao, \emph{Engineering optimization: theory and practice}.\hskip 1em plus
  0.5em minus 0.4em\relax John Wiley \& Sons, 2019.

\bibitem{kennedy1988neural}
M.~P. Kennedy and L.~O. Chua, ``Neural networks for nonlinear programming,''
  \emph{IEEE Transactions on Circuits and Systems}, vol.~35, no.~5, pp.
  554--562, 1988.

\bibitem{lillo1993solving}
W.~E. Lillo, M.~H. Loh, S.~Hui, and S.~H. Zak, ``On solving constrained
  optimization problems with neural networks: A penalty method approach,''
  \emph{IEEE Transactions on neural networks}, vol.~4, no.~6, pp. 931--940,
  1993.

\bibitem{khalil2017learning}
E.~Khalil, H.~Dai, Y.~Zhang, B.~Dilkina, and L.~Song, ``Learning combinatorial
  optimization algorithms over graphs,'' in \emph{Advances in Neural
  Information Processing Systems}, 2017, pp. 6348--6358.

\bibitem{gasse2019exact}
M.~Gasse, D.~Ch{\'e}telat, N.~Ferroni, L.~Charlin, and A.~Lodi, ``Exact
  combinatorial optimization with graph convolutional neural networks,'' in
  \emph{Advances in Neural Information Processing Systems}, 2019, pp.
  15\,554--15\,566.

\bibitem{canyasse2017supervised}
R.~Canyasse, G.~Dalal, and S.~Mannor, ``Supervised learning for optimal power
  flow as a real-time proxy,'' in \emph{2017 IEEE Power \& Energy Society
  Innovative Smart Grid Technologies Conference (ISGT)}.\hskip 1em plus 0.5em
  minus 0.4em\relax IEEE, 2017, pp. 1--5.

\bibitem{cao2018machine}
X.~Cao, R.~Ma, L.~Liu, H.~Shi, Y.~Cheng, and C.~Sun, ``A machine learning-based
  algorithm for joint scheduling and power control in wireless networks,''
  \emph{IEEE Internet of Things Journal}, vol.~5, no.~6, pp. 4308--4318, 2018.

\bibitem{cui2019spatial}
W.~Cui, K.~Shen, and W.~Yu, ``Spatial deep learning for wireless scheduling,''
  \emph{IEEE Journal on Selected Areas in Communications}, vol.~37, no.~6, pp.
  1248--1261, 2019.

\bibitem{sun2017learning}
H.~Sun, X.~Chen, Q.~Shi, M.~Hong, X.~Fu, and N.~D. Sidiropoulos, ``Learning to
  optimize: Training deep neural networks for wireless resource management,''
  in \emph{2017 IEEE 18th International Workshop on Signal Processing Advances
  in Wireless Communications (SPAWC)}.\hskip 1em plus 0.5em minus 0.4em\relax
  IEEE, 2017, pp. 1--6.

\bibitem{bojarski2016end}
M.~Bojarski, D.~Del~Testa, D.~Dworakowski, B.~Firner, B.~Flepp, P.~Goyal, L.~D.
  Jackel, M.~Monfort, U.~Muller, J.~Zhang \emph{et~al.}, ``End to end learning
  for self-driving cars,'' \emph{arXiv preprint arXiv:1604.07316}, 2016.

\bibitem{deka2019learning}
D.~Deka and S.~Misra, ``Learning for dc-opf: Classifying active sets using
  neural nets,'' in \emph{2019 IEEE Milan PowerTech}.\hskip 1em plus 0.5em
  minus 0.4em\relax IEEE, 2019, pp. 1--6.

\bibitem{ng2018statistical}
Y.~Ng, S.~Misra, L.~A. Roald, and S.~Backhaus, ``Statistical learning for dc
  optimal power flow,'' in \emph{2018 Power Systems Computation Conference
  (PSCC)}.\hskip 1em plus 0.5em minus 0.4em\relax IEEE, 2018, pp. 1--7.

\bibitem{mao2016resource}
H.~Mao, M.~Alizadeh, I.~Menache, and S.~Kandula, ``Resource management with
  deep reinforcement learning,'' in \emph{Proceedings of the 15th ACM Workshop
  on Hot Topics in Networks}, 2016, pp. 50--56.

\bibitem{song2018learning}
J.~Song, R.~Lanka, A.~Zhao, Y.~Yue, and M.~Ono, ``Learning to search via
  retrospective imitation,'' \emph{arXiv preprint arXiv:1804.00846}, 2018.

\bibitem{donti2017task}
P.~Donti, B.~Amos, and J.~Z. Kolter, ``Task-based end-to-end model learning in
  stochastic optimization,'' in \emph{Advances in Neural Information Processing
  Systems}, 2017, pp. 5484--5494.

\bibitem{djolonga2017differentiable}
J.~Djolonga and A.~Krause, ``Differentiable learning of submodular models,'' in
  \emph{Advances in Neural Information Processing Systems}, 2017, pp.
  1013--1023.

\bibitem{lee2019meta}
K.~Lee, S.~Maji, A.~Ravichandran, and S.~Soatto, ``Meta-learning with
  differentiable convex optimization,'' in \emph{Proceedings of the IEEE
  Conference on Computer Vision and Pattern Recognition}, 2019, pp.
  10\,657--10\,665.

\bibitem{srikant2013communication}
R.~Srikant and L.~Ying, \emph{Communication networks: an optimization, control,
  and stochastic networks perspective}.\hskip 1em plus 0.5em minus 0.4em\relax
  Cambridge University Press, 2013.

\bibitem{gross2005graph}
J.~L. Gross and J.~Yellen, \emph{Graph theory and its applications}.\hskip 1em
  plus 0.5em minus 0.4em\relax CRC press, 2005.

\bibitem{blumensath2008iterative}
T.~Blumensath and M.~E. Davies, ``Iterative thresholding for sparse
  approximations,'' \emph{Journal of Fourier analysis and Applications},
  vol.~14, no. 5-6, pp. 629--654, 2008.

\bibitem{norouzi2012hamming}
M.~Norouzi, D.~J. Fleet, and R.~R. Salakhutdinov, ``Hamming distance metric
  learning,'' in \emph{Advances in neural information processing systems},
  2012, pp. 1061--1069.

\bibitem{jansen1997sensitivity}
B.~Jansen, J.~De~Jong, C.~Roos, and T.~Terlaky, ``Sensitivity analysis in
  linear programming: just be careful!'' \emph{European Journal of Operational
  Research}, vol. 101, no.~1, pp. 15--28, 1997.

\bibitem{athay1979practical}
T.~Athay, R.~Podmore, and S.~Virmani, ``A practical method for the direct
  analysis of transient stability,'' \emph{IEEE Transactions on Power Apparatus
  and Systems}, no.~2, pp. 573--584, 1979.

\bibitem{Mittelmann17}
H.~D. Mittelmann, ``Latest benchmarks of optimization software,'' in
  \emph{INFORMS Annual Meeting}, 2017.

\bibitem{diamond2016cvxpy}
S.~Diamond and S.~Boyd, ``Cvxpy: A python-embedded modeling language for convex
  optimization,'' \emph{The Journal of Machine Learning Research}, vol.~17,
  no.~1, pp. 2909--2913, 2016.

\bibitem{andersen2013cvxopt}
M.~S. Andersen, J.~Dahl, and L.~Vandenberghe, ``Cvxopt: A python package for
  convex optimization, version 1.1. 6,'' \emph{Available at cvxopt. org},
  vol.~54, 2013.

\end{thebibliography}
